\documentclass[10pt]{article}
\usepackage{amssymb}   
\usepackage{latexsym}         
\usepackage{amscd}   
\usepackage{amsmath}          
\usepackage{amsfonts}
\usepackage[dvips]{color}

\newtheorem{definition}{Definition}[section]
\newtheorem{theorem}[definition]{Theorem}
\newtheorem{corollary}[definition]{Corollary}
\newtheorem{lemma}[definition]{Lemma}

\newtheorem{assumption}[definition]{Assumption}

\newtheorem{notation}[definition]{Notation}

\def\Mxc{{\mbox{Mat}}_{X}(\mathbb{C})}
\newcommand{\proof}{\noindent {\it Proof:} \ }

\begin{document}
\title{\bf Structure of Thin Irreducible Modules of a $Q$-polynomial Distance-Regular Graph}

\author{
Diana R. Cerzo{\footnote{
Institute of Mathematics,
University of the Philippines,
Diliman, Quezon City,
Philippines
}}}

\date{}

\maketitle

\begin{abstract}
Let $\Gamma$ be a $Q$-polynomial distance-regular graph with vertex set $X$, diameter $D \geq 3$ and adjacency matrix $A$.  Fix $x \in X$ and let $A^*=A^*(x)$ be the corresponding dual adjacency matrix.  Recall that the Terwilliger algebra $T=T(x)$ is the subalgebra of $\Mxc$ generated by $A$ and $A^*$.  Let $W$ denote a thin irreducible $T$-module.  It is known that the action of $A$ and $A^*$ on $W$ induces a linear algebraic object known as a Leonard pair.  Over the past decade, many results have been obtained concerning Leonard pairs.  In this paper, we apply these results to obtain a detailed description of $W$.  In our description, we do not assume that the reader is familiar with Leonard pairs.  Everything will be proved from the point of view of $\Gamma$.

Our results are summarized as follows.  Let $\{E_i\}_{i=0}^D$ be a $Q$-polynomial ordering of the primitive idempotents of $\Gamma$ and let $\{E^*_{i}\}_{i=0}^D$ be the dual primitive idempotents of $\Gamma$ with respect to $x$.  Let $r, t$ and $d$ be the endpoint, dual endpoint and diameter of $W$, respectively.  Let  $u \mbox{ and }v$ be nonzero vectors in $E_tW$ and $E^*_rW\!,$ respectively. We show that $\{E^*_{r+i}A^iv\}_{i=0}^d$ and $\{E_{t+i}A^{*i}u \}_{i=0}^d$ are bases for $W$ that are orthogonal with respect to the standard Hermitian dot product.  We display the matrix representations of $A$ and $A^*$ with respect to these bases.  We associate with $W$ two sequences of polynomials $\{p_i\}_{i=0}^d \mbox{ and } \{p^*_i\}_{i=0}^d$.  We show that  for $0 \leq i \leq d$, $p_i(A)v = E^*_{r+i}A^iv$ and $p^*_i(A^*)u = E_{t+i}A^{*i}u.$  Next, we show that $\{E^*_{r+i}u \}_{i=0}^d$ and $\{E_{t+i}v \}_{i=0}^d$ are orthogonal bases for $W\!$; we call these the standard basis and dual standard basis for $W\!$, respectively. We display the matrix representations of $A$ and $A^*$ with respect to these bases.  The entries in these matrices will play an important role in our theory.  We call these the intersection numbers and dual intersection numbers of $W\!.$  Using these numbers, we compute all inner products involving the standard and dual standard bases.   We also use these numbers to define two normalizations $u_i, v_i$ (resp.\! $u^*_i, v^*_i$) for $p_i$ (resp.\! $p^*_i$).  Using the orthogonality of the standard and dual standard bases, we show that for each of the sequences $\{p_i\}_{i=0}^d, \{p^*_i\}_{i=0}^d,  \{u_i\}_{i=0}^d,  \{u^*_i\}_{i=0}^d, \{v_i\}_{i=0}^d, \{v^*_i\}_{i=0}^d$ the polynomials involved are orthogonal and we display the orthogonality relations.  We also show that each of the sequences satisfy a three-term recurrence and a relation known as the Askey-Wilson duality.  We then turn our attention to two more bases for $W\!.$  We find the matrix representations of $A$ and $A^*$ with respect to these bases.  From the entries of these matrices we obtain two sequences of scalars known as the first split sequence and second split  sequence of $W\!.$  We associate with $W$ a sequence of scalars called the parameter array.  This sequence consists of the eigenvalues of the restriction of $A$ to $W$, the eigenvalues of the restriction of $A^*$ to $W\!$, the first split sequence of $W$ and the second split sequence of $W$. We express  all the scalars and polynomials associated with $W$ in terms of its parameter array.  We show that the parameter array of $W$ is determined by $r, t, d$ and one more free parameter.  From this we conclude that the isomorphism class of $W$ is determined by these four parameters.   Finally, we apply our results to the case in which $\Gamma$ has $q$-Racah type or classical parameters.       
  
\end{abstract}

\section{Introduction}
\indent The Terwilliger algebra $T$ of a distance-regular graph was first introduced in \cite{subI}.  This algebra has been used extensively to study the $Q$-polynomial property \cite{Caugh,Curtin,Go}.  In this paper, we continue this study focusing on the structure of thin irreducible $T$-modules.

Let $\Gamma$ be a $Q$-polynomial distance-regular graph with vertex set $X$, diameter $D \geq 3,$ and adjacency matrix $A$ (see Section 2 for formal definitions). Fix $x \in X$ and let $A^* = A^*(x)$ be the corresponding dual adjacency matrix.  Recall that the Terwilliger algebra $T=T(x)$ is the subalgebra of $\Mxc$ generated by $A$ and $A^*$.  Let $W$ be a thin irreducible $T$-module.  It is known that the action of $A$ and $A^*$ on $W$ induces a linear algebraic object called a Leonard pair; this was first introduced by Terwilliger in \cite{split}.  The theory of Leonard pairs has been developed over the past decade. We apply these results to obtain a detailed description of $W$. In our description, we do not assume that the reader is familiar with Leonard pairs.  The results will be proved from the point of view of $\Gamma$.

Our results are summarized as follows.  Let $\{E_i\}_{i=0}^D$ be a $Q$-polynomial ordering of the primitive idempotents of $\Gamma$ and let $\{E^*_{i}\}_{i=0}^D$ be the dual primitive idempotents of $
\Gamma$ with respect to $x$.  Let $r, t$ and $d$ be the endpoint, dual endpoint and diameter of $W$, respectively.  Let  $u \mbox{ and }v$ be nonzero vectors in $E_tW$ and $E^*_rW\!,$ respectively. We show that $\{E^*_{r+i}A^iv\}_{i=0}^d$ and $\{E_{t+i}A^{*i}u \}_{i=0}^d$ are bases for $W$ that are orthogonal with respect to the standard Hermitian dot product.  We display the matrix representations of $A$ and $A^*$ with respect to these bases.  We associate with $W$ two sequences of polynomials $\{p_i\}_{i=0}^d \mbox{ and } \{p^*_i\}_{i=0}^d$.  We show that  for $0 \leq i \leq d$, $p_i(A)v = E^*_{r+i}A^iv$ and $p^*_i(A^*)u = E_{t+i}A^{*i}u.$  Next, we show that $\{E^*_{r+i}u \}_{i=0}^d$ and $\{E_{t+i}v \}_{i=0}^d$ are orthogonal bases for $W\!$; we call these the standard basis and dual standard basis for $W\!$, respectively. We display the matrix representations of $A$ and $A^*$ with respect to these bases.  The entries in these matrices will play an important role in our theory.  We call these the intersection numbers and dual intersection numbers of $W\!.$  Using these numbers, we compute all inner products involving the standard and dual standard bases.   We also use these numbers to define two normalizations $u_i, v_i$ (resp.\! $u^*_i, v^*_i$) for $p_i$ (resp.\! $p^*_i$).  Using the orthogonality of the standard and dual standard bases, we show that for each of the sequences $\{p_i\}_{i=0}^d, \{p^*_i\}_{i=0}^d,  \{u_i\}_{i=0}^d,  \{u^*_i\}_{i=0}^d, \{v_i\}_{i=0}^d, \{v^*_i\}_{i=0}^d$ the polynomials involved are orthogonal and we display the orthogonality relations.  We also show that each of the sequences satisfy a three-term recurrence and a relation known as the Askey-Wilson duality.  We then turn our attention to two more bases for $W\!.$  We find the matrix representations of $A$ and $A^*$ with respect to these bases.  From the entries of these matrices we obtain two sequences of scalars known as the first split sequence and second split  sequence of $W\!.$  We associate with $W$ a sequence of scalars called the parameter array.  This sequence consists of the eigenvalues of the restriction of $A$ to $W$, the eigenvalues of the restriction of $A^*$ to $W\!$, the first split sequence of $W$ and the second split sequence of $W$. We express  all the scalars and polynomials associated with $W$ in terms of its parameter array.  We show that the parameter array of $W$ is determined by $r, t, d$ and one more free parameter.  From this we conclude that the isomorphism class of $W$ is determined by these four parameters.   Finally, we apply our results to the case in which $\Gamma$ has $q$-Racah type or classical parameters.       
  
\section {Preliminaries} In this section, we recall some basic concepts concerning $Q$-polynomial distance-regular graphs.    For more background information see 
\cite{BI} and \cite{BCN}.\\
\indent
Let $X$ be a non-empty finite set.  
Let Mat${}_X(\mathbb{C})$ denote the $\mathbb{C}$-algebra of matrices
whose rows and columns are indexed by $X$ and whose entries are in $\mathbb{C}$.  We let $I$ (resp.\! $J$) denote the identity matrix (resp.\! all 1's matrix) in $\Mxc$.  
Let $V=\mathbb{C}X$ be the vector space over $\mathbb{C}$ consisting of 
column vectors whose coordinates are indexed by $X$ and 
whose entries are in $\mathbb{C}$.  Observe that Mat${}_X(\mathbb{C})$
acts on $V$ by left multiplication.
For $u, v \in V$, define 
$\langle u, v \rangle := u^t\overline{v}$, where $u^t$ is the transpose of $u$ and $\overline{v}$ is the complex conjugate of $v$.  
Observe that $\langle \ , \ \rangle$ is a positive definite Hermitian form on $V$.  Note that $\langle Bu, v \rangle = \langle u, \overline{B}^t v \rangle$ 
for all $B \in$ Mat${}_X(\mathbb{C})$ and $u,v \in V$.  For $y \in X$, let $\hat{y}$ denote the element in $V$ with a $1$ in the $y$ coordinate and $0$ in all other coordinates.  Observe that $\{\hat{y} \mid y \in X \}$ is an orthonormal basis for $V\!.$\\
\indent
Let $\Gamma = (X,R)$ be a finite undirected connected graph without loops or multiple edges, with vertex set $X$ and edge set $R$.  Let $\partial$ denote the path-length distance function for $\Gamma$.  Set $D= \max \{\partial(x,y) \mid x, y \in X \}$.  We refer to $D$ as the \textit{diameter} of $\Gamma$.  For $x \in X$ and an integer $i \geq 0$, let $\Gamma_i(x) = \{y \mid y \in X, \ \partial(x,y) =i \}$.  Abbreviate $\Gamma(x) := \Gamma_1(x)$.  For an integer $k \geq 0$, we say that $\Gamma$ is \textit{regular} with \textit{valency} $k$ whenever $k=|\Gamma(x)|$ for all $x \in X$.    We say that $\Gamma$ is \textit{distance-regular} whenever there exists scalars $p^h_{ij}$ $(0 \leq h, i, j \leq D)$ such that  $p^{h}_{ij}=|\Gamma_i(x) \cap \Gamma_j(y)|$ for all $x, y \in X$ with $\partial (x,y) = h$.  We refer to the $p^h_{ij}$ as the \textit{intersection numbers} of $\Gamma$.  For the rest of this paper, assume that $\Gamma$ is distance-regular with diameter $D \geq 3$. Note that by the triangle inequality, we have 
(i) $p^h_{ij} = 0 $ if one of $h, i, j$ is greater than the sum of the other two; (ii) $p^h_{ij}\neq 0$ if one of $h, i, j$ is equal to the sum of the other two. 
We abbreviate $c_i: = p^i_{1i-1} \ (1 \leq i \leq D), \ a_i := p^i_{1i} \ (0 \leq i \leq D), \ b_i:= p^i_{1i+1} \ (0 \leq i \leq D-1)$.  For notational convenience, define $b_D=0, \ c_0 = 0$.  Observe that $\Gamma$ is regular with valency $k=b_0$.  To avoid trivialities, we always assume that $k  \geq 3$. Note that $c_i + a_i + b_i = k$ for $0 \leq i \leq D$.  
For $0 \leq i \leq D$, let $k_i = p^0_{ii}.$  Observe that $k_i = |\Gamma_i(x)|$ for all $x \in X$.  By \cite[p.195]{BI}, 
\begin{equation}
k_i = {\displaystyle\frac{b_0b_1 \cdots b_{i-1}}{c_1c_2 \cdots c_i}} \quad (0 \leq i \leq D). \label{valencyform}
\end{equation}
We refer to $k_i$ as the $i$\textit{th valency} of $\Gamma\!.$

We now recall the Bose-Mesner algebra of $\Gamma$.  
For $0 \leq i \leq D$, define $A_i \in \mbox{Mat}{}_{X}(\mathbb{C})$ to have $(x,y)$-entry equal to 1 if $\partial(x,y) = i$, and $0$ otherwise.
We refer to $A_i$ as the $i$\textit{th distance matrix} of $\Gamma$. Note that (i) $A_0  =  I$; 
(ii) $\sum_{i=0}^D A_i  =  J$;
(iii) $A_i^t  =   A_i \ (0 \leq i \leq D )$; 
(iv) $A_iA_j  =   \sum_{h = 0}^Dp^h_{ij}A_h \ (0 \leq i, j \leq D)$.   
Observe that $\{A_i \}_{i=0}^D$ are linearly independent.  Thus, they form a basis for a subalgebra $M$ of Mat${}_X(\mathbb{C})$; $M$ is called the \textit{Bose-Mesner algebra of $\Gamma$}.  Abbreviate $A:=A_1$ and call this the \textit{adjacency matrix} of $\Gamma$.  By \cite[p.190]{BI}, $M$ is generated by $A$.  By \cite[p.59]{BI},
$M$ has a second basis $\{E_i\}_{i=0}^D$ 
which satisfies the following: 
(i) $E_0 = |X|^{-1}J$; 
(ii) $\sum_{i=0}^D E_i =   I$; 
(iii) $E_i^t  =  E_i = \overline{E_i} \ (0 \leq i \leq D)$; 
(iv) $E_iE_j  =  \delta_{ij}E_i \ (0 \leq i, j \leq D)$. For notational convenience, define $E_{-1} =0, \ E_{D+1}= 0$. 
For $0 \leq i \leq D$, let $m_i$ denote the rank of $E_i$; we call $m_i$ the \textit{multiplicity} of $\Gamma$ associated with $E_i$.  
Since $\{E_i\}_{i=0}^D$ is a basis for $M$, there exist complex scalars 
$\{\theta_i\}_{i=0}^D$ such that 
$A  = \sum_{i=0}^D \theta_i E_i$.  
Note that for 
$0 \leq i \leq D$, $AE_i = E_iA = \theta_iE_i.$ Thus, $E_iV$ is an eigenspace for $A$, and $\theta_i$ is the corresponding eigenvalue.  Since $A$ is symmetric, 
$\theta_i \in \mathbb{R}$.  Since $A$ generates $M$, the $\{\theta_i \}_{i=0}^D$ are mutually distinct.  Note that    
\begin{equation}
V = \sum_{i=0}^D E_iV \quad \mbox{(orthogonal direct sum),} \label{eq:decomposeV}
\end{equation}
and that   
\begin{equation}
E_i = \prod_{\begin{subarray}{c}
0 \leq j \leq D\\
j \neq i
\end{subarray}}\frac{A-\theta_jI}{\theta_i - \theta_j} \qquad (0 \leq i \leq D). \label{eq:explicitEi}
\end{equation} 
We call $\theta_i$ the \textit{eigenvalue} of $\Gamma$ associated with $E_i$.  

We now recall the Krein parameters of $\Gamma$.  Observe that $A_i \circ A_j = \delta_{ij}A_i$ for $0 \leq i, j \leq D$, where $\circ$ is the entry-wise multiplication.  Thus, $M$ is closed under $\circ$.  
Consequently, there exist complex scalars $q_{ij}^h \ (0 \leq h, i, j \leq D)$  such that
$$E_i\circ E_j  =  |X|^{-1}\sum_{h = 0}^D q^h_{ij} E_h \quad (0 \leq i, j \leq D).$$ 
The $q^h_{ij}$ are known as the 
\textit{Krein parameters} or \textit{dual intersection numbers} of $\Gamma$.  By 
\cite[p.69]{BI}, the
$q^h_{ij}$ 
are real and nonnegative.       \\
\indent
We now consider the $Q$-polynomial property.  
The graph $\Gamma$ is said to be \textit{$Q$-polynomial} 
(with respect to the given ordering 
$\{E_i \}_{i=0}^D$ of primitive idempotents) whenever both: 
(i) $q^h_{ij}=0$ if one of $h, i, j$ is 
greater than the sum of the other two; 
(ii) $q^h_{ij} \neq 0$ if one of 
$h, i, j$ is equal to the sum of the other two.  For the rest of this paper, we assume that $\Gamma$ is $Q$-polynomial with respect to $\{E_i \}_{i=0}^D$. We abbreviate $c^*_i: = q^i_{1i-1} \ (1 \leq i \leq D), \ a^*_i := q^i_{1i} \ (0 \leq i \leq D), \ b^*_i:= q^i_{1i+1} \ (0 \leq i \leq D-1)$. For notational convenience, define $b^*_D=0, \ c^*_0=0$.  By \cite[p.67]{BI}, $m_i = q^0_{ii} \ (0 \leq i \leq D)$.  By \cite[p.196]{BI}, 
\begin{equation}
m_i = {\displaystyle\frac{b^*_0b^*_1 \cdots b^*_{i-1}}{c^*_1c^*_2 \cdots c^*_i}} \quad (0 \leq i \leq D). \label{multform}
\end{equation}
\indent
We now recall the dual Bose-Mesner algebra of $\Gamma$.  For the rest of this paper, fix $x \in X$.  
For $0 \leq i \leq D$, define 
$E^*_i= E^*_i(x)$ 
to be the diagonal matrix in Mat$_{X}(\mathbb{C})$ 
with $(y,y)$-entry
\begin{equation}
(E^*_i)_{yy} = \left\{\begin{array}{ll}
                1 & \mbox{if }\partial(x,y) = i\\
                0 & \mbox{otherwise} \label{defei*}
                \end{array}\right.  \qquad (y \in X).
\end{equation}
We refer to $E^*_i$ as the $i$\textit{th} \textit{dual primitive idempotent} of 
$\Gamma$ with respect to $x$.  For notational convenience, define $E^*_{-1}=0, \ E^*_{D+1}=0$.
Note that 
(i) $\sum_{i=0}^D E^*_i = I$; 
(ii) $E^{*t}_i = E^*_i =  \overline{E^*_i} \ (0 \leq i \leq D)$;
(iii)   $E_i^*E^*_j = \delta_{ij}E^*_i \ (0 \leq i, j \leq D)$.
Observe that $\{E^*_i \}_{i=0}^D$ are linearly independent.  Thus, they 
form a basis for a commutative subalgebra $M^* = M^*(x)$ of Mat${}_{X}(\mathbb{C})$; $M^*$ is called the 
\textit{dual Bose-Mesner algebra} of $\Gamma$ with respect to $x$.  
For $0 \leq i \leq D$, define 
$A^*_i = A^*_i(x) $ 
to be the diagonal matrix in Mat$_{X}(\mathbb{C})$ such that 
${(A^*_i)}_{yy} = |X|(E_i)_{xy}$ for $y \in X$.  By 
\cite[p.379]{subI}, 
$\{A^*_i\}_{i=0}^D$ is a basis for $M^*$ 
and satisfies the following properties: 
(i) $A^*_0 = I$; 
(ii) $ \sum_{i=0}^D A^*_i = |X|E^*_0$; 
(iii) $ A^{*t}_i= A^*_i = \overline{A^*_i} \ (0 \leq i \leq D)$; 
(iv) $A^*_iA^*_j =  \sum_{h=0}^D q^{h}_{ij} A^*_h \ (0 \leq i, j \leq D)$.  
We refer to 
$A^*_i$ as the $i$\textit{th} \textit{dual distance matrix} 
of $\Gamma$ with respect to $x$.  Abbreviate $A^* := A^*_1$ and call this the \textit{dual adjacency matrix} of $\Gamma$ with respect ot $x$.  By \cite[Lemma 3.11]{subI}, $M^*$ is generated by $A^*$.
Since $\{E^*_i\}_{i=0}^D$ is a basis for $M^*$, there exist complex scalars 
$\{\theta^*_i\}_{i=0}^D$ 
such that 
$A^*  = \sum_{i=0}^D \theta^*_i E^*_i$.  Note that for $0 \leq i \leq D$, $A^*E^*_i = E^*_iA^*= \theta^*_iE^*_i$.   Since $A^*$ is real, $\theta^*_i \in \mathbb{R}$.  Since $A^*$ generates $M^*\!$, the $\{\theta^*_i\}_{i=0}^D$ are mutually distinct.  Observe that
$$E^*_iV = \mbox{Span}\{\hat{y}\! \mid \! y \in X, \ \partial(x,y) = i \} \qquad(0 \leq i \leq D).$$ Moreover,
\begin{equation}
V = \sum_{i=0}^D E^*_iV  \qquad \mbox{(orthogonal direct sum)} \label{eq:decomposeV2}
\end{equation}
and 
\begin{equation}
E^*_i = \prod_{\begin{subarray}{c}
0 \leq j \leq D\\
j \neq i
\end{subarray}}\frac{A^*-\theta^*_jI}{\theta^*_i - \theta^*_j} \qquad (0 \leq i \leq D). \label{eq:formpidem}
\end{equation}
We call $\theta^*_i$ the \textit{dual eigenvalue} of 
$\Gamma$ associated with $E^*_i$.  

We now recall the Terwilliger algebra of $\Gamma$.  Let $T = T(x)$ denote the subalgebra of $\Mxc$ generated by $M$ and $M^*$. We refer to $T$ as the \textit{Terwilliger algebra of} $\Gamma$ \textit{with respect to} $x$. Observe that  $T$ is generated by $A,A^*$.  Moreover, $T$ is semi-simple. By \cite[Lemma 3.2]{subI},
\begin{align}
E^{*}_iA_hE^{*}_j & = 0 \mbox{ if and only if } p^{h}_{ij} = 0 \quad \quad (0 \leq h, i, j \leq D), \label{trineq1}\\ 
E_iA_h^{*}E_j & = 0 \mbox{ if and only if } q^{h}_{ij} = 0\quad \quad (0 \leq h,i, j \leq D). \label{trineq2}
\end{align}
It follows from (\ref{trineq1}) and (\ref{trineq2}) that 
\begin{align*}
AE^*_iV & \subseteq E^*_{i-1}V + E^*_{i}V + E^*_{i+1}V \qquad (0 \leq i \leq D),\\
A^*E_iV & \subseteq E_{i-1}V + E_{i}V + E_{i+1}V \qquad (0 \leq i \leq D).
\end{align*}
Moreover, 
\begin{align}
E^{*}_{i}A^hE^{*}_{j} \ &= \ \left\{\begin{array}{ll} 
                0, & h <  |\:i-j|\\
                \neq 0, & h = |\:i-j|
                \end{array}\right. \quad(0 \leq h, i, j\leq D),  \label{eq:triple1} \\ 
E_{i}A^{*h}E_{j} \ & = \ \left\{\begin{array}{ll} 
                0, & h <  |\:i-j|\\
                \neq 0, & h = |\:i-j|
                \end{array}\right. \quad (0 \leq h, i, j \leq D).\label{eq:triple2}
\end{align}
\begin{lemma} \label{lemma:products}
 For $0 \leq i, j, k, l \leq D$ with $i+j = |k-l|$,
$$E^*_{l}A^{i+j}E^*_{k} = \left\{\begin{array}{ll}
E^*_{l}A^iE^*_{l+i}A^jE^*_{k}, & i+j = k-l \\
E^*_{l}A^iE^*_{k+j}A^jE^*_{k}, & i+j = l-k.
\end{array}
\right. $$
\end{lemma}
\proof
In $E^*_lA^{i+j}E^*_k$, write $A^{i+j}$ as $A^iIA^j$ with $I = \sum_{m=0}^DE^*_m$.  Evaluate the result using (\ref{eq:triple1}).  \hfill $\Box$

\begin{lemma}
For $0 \leq i, j, k, l \leq d$ with $i+j = |k-l|$,
$$E_{l}(A^*)^{i+j}E_{k} = \left\{\begin{array}{ll}
E_{l}A^{*i}E_{l+i}A^{*j}E_{k}, & i+j = k-l \\
E_{l}A^{*i}E_{k+j}A^{*j}E_{k}, & i+j = l-k.
\end{array}
\right. $$
\end{lemma}

\proof
Similar to the proof of Lemma \ref{lemma:products}. \hfill$\Box$

\section{$T$-modules}

In this section, we recall some basic facts concerning the $T$-modules of $\Gamma$.

Let $W$ be a subspace of $V$.  We say that $W$ is a $T$\!-\textit{module} whenever $TW \subseteq W\!.$  Note that $V$ is a $T\!$-module.  We refer to $V$ as the \textit{standard module}.  Let $W$ and $W'$ be $T\!$-modules.  By a $T$-\textit{module isomorphism} from $W$ to $W'$, we mean a vector space isomorphism $\sigma: W \rightarrow W'$  such that $(\sigma B - B\sigma)W = 0$ for all $B \in T$.  If such a map exists, we say that $W$ and $W'$ are \textit{isomorphic} as $T\!$-modules. A $T\!$-module $W$ is said to be \textit{irreducible} whenever $W \neq 0$ and $W$ contains no $T$-modules besides $0$ and $W$.  
$W$ is said to be \textit{thin} whenever $\dim E^*_iW \leq 1$ for $0 \leq i \leq D$.
Similarly, $W$ is said to be \textit{dual thin} whenever $\dim E_iW \leq 1$ for $0 \leq i \leq D$.\\
\indent
We now recall the notion of endpoint, dual endpoint, diameter and dual diameter.  Observe that $W = \sum E^*_iW$ (orthogonal direct sum) where the sum is taken over all indices $i$ $(0 \leq i \leq D)$ such that $E^*_iW \neq 0$.  Similarly, $W = \sum E_iW$ (orthogonal direct sum) where the sum is taken over all indices $i$ $(0 \leq i \leq D)$ such that $E_iW \neq 0$.  Let $r =\min \{i \mid 0 \leq i \leq D, E^*_iW \neq 0 \}$ and $t =\min \{i \mid 0 \leq i \leq D, E_iW \neq 0 \}$.  We call $r$ and $t$ the \textit{endpoint} and \textit{dual endpoint} of $W$\!, respectively.  Let $d = |\{i \mid 0 \leq i \leq D, E^*_iW \neq 0 \}| - 1$ and $d^* = |\{i \mid 0 \leq i \leq D, E_iW \neq 0 \}| - 1$. We refer to $d$ and $d^*$ as the \textit{diameter} and \textit{dual diameter} of $W\!,$ respectively.

\begin{lemma}{\rm\cite[ Lemma 3.9]{subI}} \label{lemma:ds1}
Let $W$ be an irreducible $T$\!-module with endpoint $r$, dual endpoint $t$, diameter $d$ and dual diameter $d^*.\!$  Then (i)--(v) below hold.
\begin{enumerate}
\item[\rm(i)]  $AE^{*}_{i}W \subseteq E^{*}_{i-1}W + E^{*}_{i}W + E^{*}_{i+1}W \ (0 \leq i \leq D).$
\item[\rm(ii)] $E^*_iW \neq 0$ if and only if $r \leq i \leq r+d \ (0 \leq i \leq D)$.
\item[\rm(iii)] $E^*_iAE^*_jW \neq 0$ if $|i-j|=1 \ (0 \leq i, j \leq D)$.
\item[\rm(iv)]  $W = \sum_{i=0}^d E^*_{r+i}W$ {\rm(orthogonal direct sum)}.
\item[\rm(v)]  Suppose $W$ is thin.  Then $E_iW = E_iE^*_rW$ for $0 \leq i \leq D$.  Moreover, $W$ is dual thin and $d = d^*$. 
\end{enumerate}
\end{lemma}

\begin{lemma}{\rm\cite[Lemma 3.12]{subI}} \label{lemma:ds2}
Let $W$ be as in Lemma \ref{lemma:ds1}.  Then (i)--(v) below hold.
\begin{enumerate}
\item[\rm(i)]  $A^{*}E_{i}W \subseteq E_{i-1}W + E_{i}W + E_{i+1}W \ (0 \leq i \leq D)$.
\item[\rm(ii)] $E_iW \neq 0$ if and only if $t \leq i \leq t+d \ (0 \leq i \leq D)$.
\item[\rm(iii)] $E_iA^*E_jW \neq 0$ if $|i-j|=1 \ (0 \leq i, j \leq D)$.
\item[\rm(iv)]  $W = \sum_{i=0}^d E_{t+i}W$ {\rm(orthogonal direct sum)}.
\item[\rm(v)]  Suppose $W$ is dual thin.  Then $E^*_iW = E^*_iE_tW$ for $0 \leq i \leq D$.  Moreover, $W$ is thin and $d^* = d$. 
\end{enumerate}
\end{lemma}

\begin{lemma}{\rm\cite[Lemma 3.6]{subI}}
There exists a unique irreducible $T\!$-module of endpoint $0$, dual endpoint $0$ and diameter $D$.  Moreover, it is thin and dual thin.  We refer to this module as the trivial $T\!$-module.
\end{lemma}

\indent
For the rest of this paper, we will have the following assumption on $W\!.$

\begin{assumption} \label{W}
From now on, $W$ will denote a thin irreducible $T\!$-module with endpoint $r$\!, dual endpoint $t$ and diameter $d$. Unless otherwise stated, we assume that $d>0$. 
\end{assumption}

\section{Generators for End$(W)$}

With reference to Assumption \ref{W}, let End$(W)=\mbox{ End}{}_{\mathbb{C}}(W)$ denote the $\mathbb{C}$-algebra of all $\mathbb{C}$-linear transformations from $W$ to $W$\!.  In this section, we will look at bases and generators of End$(W)$. We begin with two lemmas whose proofs are routine and left to the reader.

\begin{lemma} \label{eirep}
For $0 \leq i \leq d$, let $w_i$ be a nonzero vector in $E^{*}_{r+i}W$\!.   Note that $\{w_i \}_{i=0}^{d}$ is a basis for $W\!.$  With respect to this basis,
\begin{enumerate}
\item[\rm(i)] the matrix representation of $E_{r+i}^{*}$ has $(i,i)$-entry $1$ and all other entries $0$ $(0 \leq i \leq d)$;
\item[\rm(ii)] the matrix representation of $A^{*}$ is $diag(\theta^{*}_r, \theta^{*}_{r+1}, \ldots, \theta^{*}_{r+d})$;
\item[\rm(iii)] the matrix representation of $A$ is tridiagonal with each entry nonzero on the superdiagonal and subdiagonal.
\end{enumerate}
\end{lemma}
\begin{lemma} \label{eirep*}
For $0 \leq i \leq d$, let $w^*_i$ be a nonzero vector in $E_{t+i}W\!.$   Note that $\{w^*_i \}_{i=0}^{d}$ is a basis for $W\!.$  With respect to this basis,
\begin{enumerate}
\item[\rm(i)] the matrix representation of $E_{t+i}$ has $(i,i)$-entry $1$ and all other entries $0$ $(0 \leq i \leq d)$;
\item[\rm(ii)] the matrix representation of $A$ is $diag(\theta_t, \theta_{t+1}, \ldots, \theta_{t+d})$;
\item[\rm(iii)] the matrix representation of $A^*$ is tridiagonal with each entry nonzero on the superdiagonal and subdiagonal.
\end{enumerate}
\end{lemma}

\begin{definition}
\rm We refer to the sequence $\{\theta_{t+i} \}_{i=0}^d$ (resp.\! $\{\theta^*_{r+i} \}_{i=0}^d$) as the \textit{eigenvalue sequence} (resp.\! \textit{dual eigenvalue sequence}) of $W\!.$
\end{definition}

\begin{lemma} \label{0}
On $W\!$,
$$\prod_{i=0}^d(A-\theta_{t+i}I)=0, \qquad  \qquad \prod_{i=0}^d(A^*-\theta^*_{r+i}I)=0.$$
\end{lemma}

\proof
Immediate from Lemmas \ref{eirep}(ii) and \ref{eirep*}(ii). \hfill $\Box$

\begin{lemma} \label{entries}
Let $B$ (resp.\! $B^*$) denote the matrix representation of $A$ (resp.\! $A^*$) with respect to the basis given in Lemma \ref{eirep} (resp.\! Lemma \ref{eirep*}).  Then
\begin{align*}
{(B^{h})}_{ij} \ &= \ \left\{\begin{array}{ll}
                0, & h <  |i-j|\\
                \neq 0, & h = |i-j|
                \end{array}\right. \quad (0 \leq h,i,j \leq d),\\
{(B^{*h})}_{ij} \ &= \ \left\{\begin{array}{ll}
                0, & h <  |i-j|\\
                \neq 0, & h = |i-j|
                \end{array}\right. \quad (0 \leq h,i,j \leq d).
\end{align*}
\end{lemma}

\proof
Routine using Lemmas \ref{eirep}(iii) and \ref{eirep*}(iii). \hfill $\Box$\\

Using Lemma \ref{entries} we obtain the following strengthening of (\ref{eq:triple1}) and (\ref{eq:triple2}).

\begin{lemma}\label{relations}
For $0 \leq h,i,j \leq d$, the following hold on $W\!.$
\begin{align}
E^{*}_{r+i}A^hE^{*}_{r+j} \ &= \ \left\{\begin{array}{ll}
                0, & h <  |i-j|\\
                \neq 0, & h = |i-j|,
                \end{array}\right.  \label{triple1}\\
E_{t+i}A^{*h}E_{t+j} \ &= \ \left\{\begin{array}{ll}
                0, & h <  |i-j|\\
                \neq 0, & h = |i-j|.
               \end{array}\right. \label{triple2}
\end{align}
\end{lemma}

\proof
Let $B$ denote matrix representation of $A$ with respect to the basis given in Lemma \ref{eirep}. By construction, the matrix representation of  $E^{*}_{r+i}A^hE^{*}_{r+j}$ with respect to this basis has $(i,j)$-entry $(B^h)_{ij}$ and all other entries are $0$.  Line (\ref{triple1}) follows from this and Lemma \ref{entries}.  The proof of (\ref{triple2}) is similar.  \hfill$\Box$   

\begin{theorem} \label{theorembasis}
Each of the following forms a basis for the $\mathbb{C}$-vector space End$(W)$:
\begin{enumerate}
\item[\rm(i)]  the actions of $\{A^{m}E_{r}^{*}A^n \mid 0 \leq m,  n \leq d  \}$ on $W\!$,
\item[\rm(ii)] the actions of $\{A^{*m}E_{t}A^{*n}\mid 0 \leq m,  n \leq d \}$ on $W$.
\end{enumerate}
\end{theorem}

\proof
Let $S$ denote $\{A^{m}E_{r}^{*}A^n \mid 0 \leq m,  n \leq d  \}$. Observe that $|S| =(d+1)^2$ and this is equal to the dimension of End$(W)$.  It suffices to show that the actions of the elements of $S$ on $W$ are linearly independent.  Let $\{w_i \}_{i=0}^d$ be the basis for $W$ in Lemma \ref{eirep}.  With respect to this basis, let $B$ and $F^*_r$ be the matrix representations of $A$ and $E^*_r.$  We claim that for $0 \leq m, n \leq d$, $B^mF_r^*B^n$ has entries
\begin{equation}
(B^{m}F_r^{*}B^{n})_{ij} =  \left\{\begin{array}{ll}
                 0, & i > m \mbox{ or } j > n\\
                \neq 0, & i = m \mbox{ and } j = n
                \end{array}\right. \quad (0 \leq  i, j \leq d). \label{eq:claim} 
\end{equation}
By Lemma \ref{eirep}(i), $F^{*}_r$ has $(0,0)$-entry 1 and all other entries are $0$.  Thus,            $$(B^{m}F_r^{*}B^{n})_{ij} = {(B^m)}_{i0}{(B^n)}_{0j} \qquad (0 \leq i,j \leq d).$$ 
Combining this with Lemma \ref{entries}, we obtain (\ref{eq:claim}).  It follows from (\ref{eq:claim}) that actions of the elements of $S$ on $W$ are linearly independent and hence form a basis for End$(W)$.  Similarly, (ii) can be shown to be a basis for $\mbox{End}(W)$. \hfill $\Box$

\begin{theorem}\label{generators}
Each of the following is a generating set for the $\mathbb{C}$-algebra End$(W)$:
\begin{enumerate}
\item[\rm(i)] the actions of $A, E^*_r$ on $W\!,$
\item[\rm(ii)] the actions of  $A^*\!,  E_t$ on $W\!,$
\item[\rm(iii)] the actions of $A, A^*$ on $W\!.$  
\end{enumerate}
\end{theorem}

\proof
By Theorem \ref{theorembasis}, (i) and (ii) are generating sets for End$(W)$. The set (iii) is a generating set for End$(W)$ by (i) and since $E^*_r$ is a polynomial in $A^*$. \hfill $\Box$

\begin{definition}\label{def:D} \rm
Define $\mathcal{D}$ (resp.\! $\mathcal{D}^*$) to be the subalgebra of End$(W)$ generated by the action of $A$ (resp. $A^*$) on $W\!$.
\end{definition}

\begin{lemma} \label{basisforD}
Each of the following forms a basis for the $\mathbb{C}$-vector space $\mathcal{D}$:
\begin{enumerate}
\item[\rm(i)] the actions of $\{A^i \}_{i=0}^d$ on $W\!,$
\item[\rm(ii)] the actions of $\{E_{t+i}\}_{i=0}^d$ on $W\!.$ 
\end{enumerate}
\end{lemma}

\proof
(i) By Lemma \ref{entries}, $\{A^i \}_{i=0}^d$ are linearly independent on $W$.  Combining this with Lemma \ref{0}, we obtain the result.\\
(ii)  Immediate from (\ref{eq:explicitEi}) and (i).  \hfill $\Box$

\begin{lemma} \label{basisforD*}
Each of the following forms a basis for the $\mathbb{C}$-vector space $\mathcal{D^*}$:
\begin{enumerate}
\item[\rm(i)] the actions of $\{A^{*i} \}_{i=0}^d$ on $W\!,$
\item[\rm(ii)] the actions of $\{E^*_{r+i}\}_{i=0}^d$ on $W\!.$ 
\end{enumerate}
\end{lemma}

\proof
Similar to the proof of Lemma \ref{basisforD}. \hfill $\Box$

\begin{corollary} \label{basisforendw}
Each of the following forms a basis for the $\mathbb{C}$-vector space End$(W)$:
\begin{enumerate}
\item[\rm(i)] the actions of $\{E_{t+i}E^{*}_rE_{t+j} \mid 0 \leq i,j \leq d\}$ on $W\!,$
\item[\rm(ii)] the actions of $\{E^*_{r+i}E_tE^*_{r+j} \mid 0 \leq i, j \leq d\}$ on $W\!.$ 
\end{enumerate}
\end{corollary}

\proof Immediate from Theorem \ref{theorembasis} and Lemmas \ref{basisforD}, \ref{basisforD*}. \hfill $\Box$

\section {The scalars $a_i(W)$ and $x_i(W)$}
Let $W$ be as in Assumption \ref{W}.  In this section, we associate with $W$ two sequences of scalars called the $a_i(W)$ and $x_i(W)$.  We will then describe the algebraic properties of these scalars.

\begin{notation} \rm
For any $Y \in T$,  tr${}_{W}Y$ denotes the trace of the action of $Y$ on $W\!.$
\end{notation}

\begin{definition} \rm \label{defxiai} 
Define
\begin{align}
a_i(W) \ &= \ \mbox{ \rm tr}_{W}(E^{*}_{r+i}A) \qquad \qquad \quad \quad  a^{*}_i(W) \ = \ \mbox{ \rm tr}_{W}(E_{t+i}A^{*})  \qquad \qquad \quad \qquad (0 \leq i \leq d),\label{defai}\\
x_i(W) \ &= \ \mbox{ \rm tr}_{W}(E^{*}_{r+i}AE^{*}_{r+i-1}A)  \qquad x^{*}_i(W)\ = \ \mbox{ \rm tr}_{W}(E_{t+i}A^{*}E_{t+i-1}A^*) \quad \qquad(1 \leq i \leq d).  \label{xi}
\end{align}
For notational convenience, define $x_0(W)=0$ and $x^*_0(W)=0. $
\end{definition}

\begin{lemma} \label{genbasis}
For $0 \leq i \leq d$, let $w_i$ be a nonzero vector in $E^{*}_{r+i}W$\!. Let $B$ denote the matrix representation of $A$ with respect to $\{w_i \}_{i=0}^d$.  Then (i)--(iii) below hold.
\begin{enumerate}
\item[\rm(i)] $B_{ii} = a_i(W) \quad (0 \leq i \leq d)$.
\item[\rm(ii)] $B_{i,i-1}B_{i-1,i} = x_i(W) \quad (1 \leq i \leq d)$.
\item[\rm(iii)] $x_i(W) \neq 0 \quad (1 \leq i \leq d)$.
\end{enumerate}
\end{lemma}

\proof (i) By Lemma \ref{eirep}(i), (iii), the $(j,j)$-entry of the matrix representation of $E^*_{r+i}A$ with respect to $\{w_i \}_{i=0}^d$ is $B_{ii}$ if $j=i$ and $0$ otherwise ($0 \leq j \leq d$).  Taking the trace of this matrix and using (\ref{defai}), we obtain the desired result.\\
(ii) By Lemma \ref{eirep}(i), (iii), the $(j,j)$-entry of the matrix representation of $E^{*}_{r+i}AE^{*}_{r+i-1}A$ with respect to $\{w_i \}_{i=0}^d$ is $B_{i, i-1}B_{i-1,i}$ if $j=i$ and $0$ otherwise ($0 \leq j \leq d$).  Taking the trace of this matrix and using (\ref{xi}), we obtain the desired result.  \\
(iii) Immediate from (ii) and Lemma \ref{eirep}(iii). \hfill $\Box$

\begin{lemma} \label{genbasis*}
For $0 \leq i \leq d$, let $w^*_i$ be a nonzero vector in $E_{t+i}W$\!. Let $B^{*}$ denote the matrix representation of $A^*$ with respect to $\{w^*_i \}_{i=0}^d$.  Then (i)--(iii) below hold.
\begin{enumerate}
\item[\rm(i)] $B^{*}_{ii} = a^{*}_i(W)\quad (0 \leq i \leq d)$.
\item[\rm(ii)] $B^{*}_{i,i-1}B^{*}_{i-1,i} = x^{*}_i(W) \quad (1 \leq i \leq d)$.
\item[\rm(iii)] $x^{*}_i(W) \neq 0 \quad (1 \leq i \leq d)$.
\end{enumerate}
\end{lemma}

\proof Similar to the proof of Lemma \ref{genbasis}.\hfill $\Box$

\begin{theorem} \label{basis}
Let $v$ be a nonzero vector in $E^{*}_rW\!$. Then for $0\! \leq i\! \leq d$, $E^{*}_{r+i}A^{i}v$ is nonzero and hence is a basis for $E^{*}_{r+i}W$\!.  Moreover,  $\{ E^{*}_{r+i}A^i v \}_{i=0}^{d}$ is a basis for $W$\!.
\end{theorem}

\proof Since $v$ spans $E^*_rW$, $E^*_{r+i}A^iv$ spans $E^*_{r+i}A^iE^*_rW$.  By Lemma \ref{relations}, $E_{r+i}^{*}A^iE_r^{*}W \neq 0$.  Hence, $E^*_{r+i}A^iv \neq 0$. The rest of the assertion follows.  \hfill $\Box$

\begin{theorem} \label{basis*}
Let $u$ be a nonzero vector in $E_tW$\!.  Then for $0 \!\leq i \leq \!d$, $E_{t+i}A^{*i}u$ is nonzero and hence is a basis for $E_{t+i}W$\!.  Moreover,  $\{E_{t+i}A^{*i}u \}_{i=0}^{d}$ is a basis for $W$\!.
\end{theorem}

\proof
Similar to the proof of Theorem \ref{basis}. \hfill $\Box$

\begin{theorem} \label{matrixrep}
With respect to the basis given in Theorem \ref{basis}, the matrix representation of $A$  is  \begin{equation}\left(
    \begin{array}{cccccc}
    a_0(W) &   x_1(W)     &          &        &           &  {\bf 0} \\
    1      &   a_1(W)     & x_2(W)   &        &           &          \\    \label{eq:matrepai}
             &   1        & \cdot    &  \cdot &           &          \\
             &            & \cdot    &  \cdot &  \cdot    &          \\  
             &            &          &  \cdot & a_{d-1}(W)& x_d(W) \\
    {\bf 0}  &            &          &        & 1         & a_d(W)  \\
    \end{array}
\right).
\end{equation}
\end{theorem}

\proof
Let $\{w_i \}_{i=0}^d$ be the basis for $W$ in Theorem \ref{basis}.  Let $B$ denote the matrix representation of $A$ with respect to this basis.  Note that for $0 \leq i \leq d-1$, $E^*_{r+i+1}Aw_i = B_{i+1,i}w_{i+1}$. By Lemma \ref{lemma:products}, $$E^*_{r+i+1}Aw_i =E^*_{r+i+1}AE^*_{r+i}A^iE^*_rv = E^*_{r+i+1}A^{i+1}E^*_rv= w_{i+1}.$$ Thus, $B_{i+1,i} = 1$ for $0 \leq i \leq d-1$.  The rest of the assertion follows from Lemma \ref{genbasis}(i), (ii). \hfill $\Box$

\begin{theorem}
With respect to the basis given in Theorem \ref{basis*}, the matrix representation of $A^{*}$  is 
$$  \left(
    \begin{array}{cccccc}
    a^{*}_0(W) &   x^{*}_1(W)     &          &        &           &  {\bf 0} \\
    1          &   a^{*}_1(W)     & x^{*}_2(W)   &        &           &          \\
             &   1        & \cdot    &  \cdot &           &          \\
             &            & \cdot    &  \cdot &  \cdot    &          \\  
             &            &          &  \cdot & a^{*}_{d-1}(W)     & x^{*}_d(W) \\
    {\bf 0}  &            &          &        & 1          & a^{*}_d(W)  \\
    \end{array}
\right).
$$
\end{theorem}
\proof
Similar to the proof of Theorem \ref{matrixrep}.  \hfill $\Box$

\begin{lemma}\label{aixiprop}
The following hold on $W$.
\begin{enumerate}
\item[\rm(i)] $E^{*}_{r+i}AE^{*}_{r+i} = a_i(W)E^{*}_{r+i} \quad (0 \leq i \leq d)$.
\item[\rm(ii)] $E^{*}_{r+i}AE^{*}_{r+i-1}AE^{*}_{r+i} = x_i(W)E^{*}_{r+i} \quad (1 \leq i \leq d)$.
\item[\rm(iii)] $E^{*}_{r+i-1}AE^{*}_{r+i}AE^{*}_{r+i-1} = x_i(W)E^{*}_{r+i-1} \quad (1 \leq i \leq d)$.
\item[\rm(iv)] $E_{t+i}A^{*}E_{t+i} = a^{*}_i(W)E_{t+i} \quad (0 \leq i \leq d)$.
\item[\rm(v)] $E_{t+i}A^{*}E_{t+i-1}A^{*}E_{t+i} = x^{*}_i(W)E_{t+i} \quad (1 \leq i \leq d)$.
\item[\rm(vi)] $E_{t+i-1}A^{*}E_{t+i}A^{*}E_{t+i-1} = x^{*}_i(W)E_{t+i-1} \quad (1 \leq i \leq d)$.
\end{enumerate}
\end{lemma}

\proof (i)  Let $\{w_j\}_{j=0}^d$ be the basis for $W$ in Theorem \ref{basis}.  By (\ref{eq:matrepai}), 
$E^*_{r+i}AE^*_{r+i}w_j = \delta_{ij}a_i(W)E^{*}_{r+i}w_j$.  The result follows.   \\
(ii) Let $G_i$ denote the action of $E^*_{r+i}$ on $W$.  Since $W$ is thin, $G_i\mbox{End}(W)G_i$ has dimension $1$.  Observe that $G_i$ is a nonzero element of $G_i\mbox{End}(W)G_i$.  Thus there exists $\alpha \in \mathbb{C}$ such that $E^{*}_{r+i}AE^{*}_{r+i-1}AE^{*}_{r+i} = \alpha E^{*}_{r+i}$ on $W$. Take the trace of both sides of this equation. Evaluating this using Definition \ref{defxiai} and the fact that tr${}_{W}(E^*_{r+i})=1$, we find that $\alpha = x_i(W)$.\\
(iii) Similar to (ii).\\
(iv)-(vi) Similar to the proofs of (i)--(iii). \hfill $\Box$

\begin{lemma} \label{sumofai}
The following hold.
\begin{enumerate}
\item[\rm(i)] $\sum_{i=0}^d a_i(W) =  \sum_{i=0}^d \theta_{t+i}$.
\item[\rm(ii)] $\sum_{i=0}^d a^*_i(W) =  \sum_{i=0}^d \theta^*_{r+i}$.
\item[\rm(iii)]  $a_i(W) \in \mathbb{R}, \ a^*_i(W) \in \mathbb{R}    \quad (0 \leq i \leq d)$.
\item[\rm(iv)]  $x_i(W) \in \mathbb{R}, \ x_i(W) > 0  \quad (1 \leq i \leq d)$.
\item[\rm(v)]  $x^*_i(W) \in \mathbb{R}, \ x^*_i(W) > 0  \quad (1 \leq i \leq d)$. 
\end{enumerate}
\end{lemma}
\proof (i)  Immediate from Theorem \ref{matrixrep} and the fact that $\{\theta_{t+i} \}_{i=0}^d$ are the eigenvalues of the action of $A$ on $W\!$.\\ 
(ii)  Similar to the proof of (i).\\
(iii) By Lemma \ref{aixiprop}(i), $a_i(W)$ is an eigenvalue of the real symmetric matrix $E^*_{r+i}AE^*_{r+i}$.  Thus, $a_i(W) \in \mathbb{R}$.  Similarly, $a^*_i(W) \in \mathbb{R}$.  \\
(iv) By Lemma \ref{aixiprop}(ii), $x_i(W)$ is an eigenvalue of the real symmetric matrix $E^{*}_{r+i}AE^{*}_{r+i-1}AE^{*}_{r+i}$.  Thus, $x_i(W) \in \mathbb{R}$.   Since $$E^{*}_{r+i}AE^{*}_{r+i-1}AE^{*}_{r+i} = (E_{r+i-1}^{*}AE_{r+i}^{*})^{t}(E_{r+i-1}^{*}AE_{r+i}^{*})$$
is positive definite, $x_i(W) > 0$.   \\
(v) Similar to the proof of (iv). \hfill $\Box$

\section{The polynomial $p_i$}

Let $W$ be as in Assumption \ref{W}.  In the previous section, we defined two bases for $W\!$.  In this section, we will use these bases to obtain two sequences of polynomials.  We will investigate some properties of these polynomials.  Let $\mathbb{C}[\lambda]$ denote the $\mathbb{C}$-algebra of polynomials in $\lambda$ with coefficients in $\mathbb{C}$.    

\begin{definition} \label{polypi}
\rm
For $0 \leq i \leq d+1$, define  $p_i = p_i^W$ in $\mathbb{C}[\lambda]$ by $p_0 =  1$,
\begin{equation}
\lambda p_{i} =  p_{i+1} + a_{i}(W)p_{i} + x_{i}(W)p_{i-1} \qquad (0 \leq i \leq d), \label{pi}
\end{equation}
where $x_i(W)$ and $a_i(W)$ are as in Definition \ref{defxiai} and $p_{-1}=0$.
\end{definition}
\begin{definition} \label{polypi*}
\rm 
For $0 \leq i \leq d+1$, define  $p^*_i = p^{*W}_i$ in $\mathbb{C}[\lambda]$ by $p^*_0  =  1$,
\begin{equation}
\lambda p^*_{i} =  p^*_{i+1} + a^*_{i}(W)p^*_{i} + x^*_{i}(W)p^*_{i-1} \qquad (0 \leq i \leq d), \label{pi*}
\end{equation}
where $x^*_i(W)$ and $a^*_i(W)$ are as in Definition \ref{defxiai} and $p^*_{-1}=0$.
\end{definition}

\begin{lemma} \label{ppush} 
For any nonzero $u \in E_tW$ and nonzero $v \in E^{*}_rW$\!,
\begin{align}
p_i(A)v \ \ &= \ \ E_{r+i}^{*}A^iv  \qquad (0 \leq i \leq d ), \label{ppush1} \\
p^{*}_i(A^*)u \ \ &= \ \ E_{t+i}A^{*i}u   \qquad (0 \leq i \leq d ). \label{ppush2}
\end{align}
Moreover, $p_{d+1}(A)v = 0$ and $p^*_{d+1}(A^*)u = 0$.
\end{lemma}

\proof For $0 \leq i \leq d+1$, let $w_i =E^{*}_{r+i}A^{i}v$ and $w_i' = p_i(A)v$.  Recall that by Theorem \ref{basis}, $\{w_i \}_{i=0}^d$ is a basis for $W\!.$  By (\ref{eq:matrepai}), 
\begin{equation}
Aw_i = w_{i+1} + a_i(W)w_i + x_i(W)w_{i-1}  \quad (0 \leq i \leq d).\label{pang2}
\end{equation}
By (\ref{pi}),
\begin{equation}
Aw'_i = w'_{i+1} + a_i(W)w'_i + x_i(W)w'_{i-1}  \quad (0 \leq i \leq d).\label{una}
\end{equation}
Comparing (\ref{pang2}) and (\ref{una}) and using the fact that $w_0 = w_0'$, we find that $w_i = w'_{i}$  for $0 \leq i \leq d+1$.  Hence $(\ref{ppush1})$ holds.  Since $w_{d+1} = 0,$ $p_{d+1}(A)v =0$.  The rest of the assertion is proved similary.   \hfill $\Box$

\begin{theorem}
For $0 \leq i \leq d$,
\begin{align}
p_i(A)E^{*}_rW \ &= \ E^{*}_{r+i}W\!, \label{first}\\
p^*_i(A^*)E_tW \ &= \ E_{t+i}W. \label{second}
\end{align}
\end{theorem}

\proof Let $v$ be a nonzero vector in $E_r^{*}W\!.$  By (\ref{ppush1}),  $E_{r+i}^{*}A^iv$ spans $p_i(A)E_r^{*}W\!.$  By Theorem \ref{basis}, $E_{r+i}^{*}A^iv$ spans $E_{r+i}^{*}W$\!.  From these comments, we obtain (\ref{first}).  The proof for (\ref{second}) is similar. \hfill $\Box$

\begin{theorem}
For $0 \leq i \leq d$, the following hold on $W$\!.
\begin{align*}
p_i(A)E^{*}_{r} \ &= \ E^{*}_{r+i}A^{i}E^{*}_r,\\
p^{*}_i(A^*)E_{t} \ &= \ E_{t+i}A^{*i}E_t.
\end{align*}
\end{theorem}

\proof Abbreviate $\Delta :=p_i(A)E^{*}_{r}- E^{*}_{r+i}A^{i}E^{*}_r$.  We will show that $\Delta = 0$ on $W$.  For $0 \leq j \leq d$, let $w_j$ be a nonzero vector in $E^*_{r+j}W\!.$  Note that $\Delta w_j = 0$ for $1 \leq j \leq d$.  By Lemma \ref{ppush}, $\Delta w_0 = 0$.  Therefore, $\Delta = 0$ on $W$.  The second assertion is proved similary. \hfill $\Box$

\begin{theorem} \label{proppi}
The following hold.
\begin{enumerate}
\item[\rm(i)] $p_{d+1}$ is both the minimal polynomial and the characteristic polynomial of the action of $A$ on $W$.
\item[\rm(ii)] $p_{d+1} = \prod_{i=0}^d (\lambda - \theta_{t+i}).$
\item[\rm(iii)] $p^*_{d+1}$ is both the minimal polynomial and characteristic polynomial of the action of $A^*$ on $W\!.$
\item[\rm(iv)] $p^*_{d+1} = \prod_{i=0}^d (\lambda - \theta^{*}_{r+i})$.
\end{enumerate}
\end{theorem}

\proof (i)  By Lemma \ref{ppush}, $p_{d+1}(A)E^*_rW = 0$.  For $1 \leq i \leq d$,   
\begin{eqnarray*}
p_{d+1}(A)E^*_{r+i}W & = & p_{d+1}(A)p_i(A)E^*_rW \qquad \mbox{by (\ref{first})} \\
& = & p_i(A) p_{d+1}(A)E^*_rW \\
& = & 0.
\end{eqnarray*}
Therefore, $p_{d+1}(A)E^*_{r+i}W=0$ for $0 \leq i \leq d$.  Hence by Lemma \ref{lemma:ds1}(iv), $p_{d+1}(A)=0$ on $W$. 
By Theorem \ref{basis} and (\ref{ppush1}), $p_i(A) \neq 0$ on $W$ for $0 \leq i \leq d$.  From these comments, $p_{d+1}$ is the minimal polynomial of the action of $A$ on $W\!$.  Since the characteristic polynomial of the action of $A$ on $W$ has degree $d+1$, it follows that $p_{d+1}$ is also the characteristic polynomial of this action.\\
(ii) Immediate from (i) and the fact that $\{\theta_{t+i}\}_{i=0}^d$ are the eigenvalues of the action of $A$ on $W\!$.\\
(iii), (iv) Similar to the proofs of (i), (ii). \hfill $\Box$

\section{The scalars $\nu, m_i$ }

Let $W$ be as in Assumption \ref{W}.  In this section, we will investigate the algebraic properties of two more scalars associated with $W,$ called the $m_i(W)$ and the $\nu(W)$. 

\begin{definition} \rm \label{defmi} 
For $0 \leq i \leq d$, define
\begin{align}
m_i(W) \ & = \ \mbox{ \rm tr}_{W}(E_{t+i}E^{*}_{r}), \label{mi}\\
m_i^{*}(W) \ & = \ \mbox{ \rm tr}_{W}(E^{*}_{r+i}E_{t}) \label{mi*}.
\end{align}
\end{definition}

\begin{lemma} \label{miprop}
For $0 \leq i \leq d$, the following (i)--(iv) hold on $W$\!.
\begin{enumerate}
\item[\rm(i)] $E_{t+i}E^{*}_rE_{t+i} = m_i(W)E_{t+i}$. 
\item[\rm(ii)] $E^{*}_rE_{t+i}E^{*}_r = m_i(W)E^{*}_r$. 
\item[\rm(iii)] $E^{*}_{r+i}E_tE^{*}_{r+i} = m^{*}_i(W)E^{*}_{r+i}$. 
\item[\rm(iv)] $E_tE^{*}_{r+i}E_t = m^{*}_i(W)E_t$. 
\end{enumerate}
\end{lemma}

\proof (i) Let $H_i$ denote the action of $E_{t+i}$ on $W$\!.  Since $W$ is thin, $H_{i}$End$(W)H_{i}$ has dimension $1$.  Note that $H_i$ is a nonzero element of $H_{i}$End$(W)H_{i}$, hence a basis for $H_{i}$End$(W)H_{i}$. Thus there exists $\alpha \in \mathbb{C}$ such that $E_{t+i}E^{*}_rE_{t+i} =  \alpha E_{t+i}$ on $W$.  Taking the trace of both sides of this equation and using Definition \ref{defmi} and the fact that tr${}_W(E^*_r) = 1$, we find that $\alpha = m_i(W)$. \\
(ii) Let $L_r$ denote the action of $E^*_r$ on $W$\!.  Since $W$ is thin, $L_r\mbox{End}(W)L_r$ has dimension $1$.  Note that $L_r$ is a nonzero element of $L_r\mbox{End}(W)L_r$, hence a basis for $L_r\mbox{End}(W)L_r$.  Thus there exists $\alpha \in \mathbb{C}$ such that $E^*_rE_{t+i}E^*_r = \alpha E^*_r$ on $W$.  Arguing as in the proof of (i), we find that $\alpha = m_i(W)$.  \\
(iii), (iv) Similar to the proofs of (i), (ii). \hfill $\Box$

\begin{lemma} \label{miprop2}
The following hold.
\begin{enumerate}
\item[\rm(i)] $\sum_{i=0}^{d} m_i(W)= 1$.
\item[\rm(ii)] $\sum_{i=0}^{d} m^{*}_i(W)=1$.
\item[\rm(iii)] $m_i(W) \in \mathbb{R}, \ m_i(W) > 0 \quad (0 \leq i \leq d)$.
\item[\rm(iv)] $m^*_i(W) \in \mathbb{R}, \ m^*_i(W) > 0 \quad (0 \leq i \leq d)$.
\end{enumerate}
\end{lemma}

\proof (i) Observe that on $W\!$, $\sum_{i=0}^d E_{t+i}=I$.  In this equation, multiply each term on the right by $E^{*}_{r}$, take the trace and use Definition \ref{defmi} to obtain $\sum_{i=0}^{d} m_i(W)= 1$.  \\
(ii) Similar to the proof of (i).\\
(iii) By Lemma \ref{miprop}(i), $m_i(W)$ is an eigenvalue of the real symmetric matrix $E_{t+i}E^{*}_rE_{t+i}$.  Hence $m_i(W) \in \mathbb{R}$.  Since $E_{t+i}E^{*}_rE_{t+i}=(E^{*}_rE_{t+i})^{t}(E^{*}_rE_{t+i})$ is positive definite, $m_i(W) > 0$. \\
(iv)  Similar to the proof of (iii). \hfill $\Box$

\begin{definition} \rm\label{defnu}
Note that $m_0(W) = m_0^*(W)$.  We denote the multiplicative inverse of this common value to be $\nu(W)$.  
\end{definition}

The following is an immediate consequence of Lemma \ref{miprop} and Definition \ref{defnu}. 
\begin{lemma}\label{thm:nu}
The following hold on $W$.
\begin{enumerate}
\item[\rm(i)] $\nu(W)E_tE^{*}_rE_t = E_{t}$. 
\item[\rm(ii)] $\nu(W)E^{*}_rE_{t}E^{*}_r = E^{*}_r$. 
\end{enumerate}
\end{lemma}

\section{Two bases for $W$} \label{sec:2bases}

Let $W$ be as in Assumption \ref{W}.  In this section, we will look at two bases for $W$ called the standard basis and dual standard basis. 

\begin{theorem}
\label{standardanddual}  Let $u$ and $v$ be nonzero vectors in $E_tW$ and $E^{*}_{r}W\!$, respectively.  Then (i)--(ii) below hold.
  \begin{enumerate}
  \item[\rm(i)] $\{E^{*}_{r+i}u \}_{i = 0}^{d}$ is a basis for $W$\!.
  \item[\rm(ii)] $\{E_{t+i}v \}_{i = 0}^{d}$ is a basis for $W$\!.
  \end{enumerate}
\end{theorem}

\proof (i) By Lemma \ref{lemma:ds1}(iv), it suffices to show that $E_{r+i}^{*}u \neq 0$ for $0 \leq i \leq d$. By Lemmas \ref{lemma:ds1}(ii) and \ref{lemma:ds2}(v), $E_{r+i}^{*}E_tW = E^*_{r+i}W \neq 0$.  Since $u$ spans $E_tW\!$, $E_{r+i}^{*}u$ spans $E_{r+i}^{*}E_tW$.  Therefore, $E_{r+i}^{*}u \neq 0$.\\
(ii) Similar to (i).\hfill $\Box$

\begin{definition} \rm \label{standardbasis}
Let $u$ and $v$ be nonzero vectors in $E_tW$ and $E^{*}_{r}W\!$, respectively.  We call $\{E^{*}_{r+i}u \}^{d}_{i=0}$ (resp. $\{E_{t+i}v \}^{d}_{i=0}$ ) a \textit{standard} (resp. \textit{dual standard}) basis for $W$.
\end{definition}

\begin{theorem} \label{thm:unique}
Let $\{w_i\}_{i=0}^{d}$ be a standard basis for $W$ and $\{w'_i\}_{i=0}^{d}$ be a sequence of vectors in $W\!.$ Then the following are equivalent.
\begin{enumerate}
\item[\rm(i)]  $\{w'_i \}^{d}_{i=0}$ is a standard basis for $W$\!.
\item[\rm(ii)] There exists  a nonzero $\alpha \in \mathbb{C}$ such that $w'_i = \alpha w_i$ for $0 \leq i \leq d$.
\end{enumerate}
\end{theorem}

\proof By Definition \ref{standardbasis}, there exists a nonzero $u \in E_tW$ such that $w_i = E^*_{r+i}u $   for $0 \leq i \leq d$.  Note that $\{w_i' \}$ is a standard basis for $W$ if and only if there exists a nonzero $u' \in E_tW$ such that $w_i' = E_{r+i}^{*}u'$ for $0 \leq  i  \leq d$.  Since $u$ spans $E_tW$\!, $u' = \alpha u$ for some nonzero $\alpha \in \mathbb{C}$.  The conclusion follows. \hfill $\Box$

\begin{theorem} \label{thm:unique*}
Let $\{v_i\}^{d}_{i=0}$ be a dual standard basis for $W$ and $\{v'_i \}^{d}_{i=0}$ be a sequence of vectors in $W$\!. Then the following are equivalent.
\begin{enumerate}
\item[\rm(i)]  $\{v'_i \}^{d}_{i=0}$ is a dual standard basis for $W$\!.
\item[\rm(ii)] There exists a nonzero $\alpha \in \mathbb{C}$ such that $v'_i = \alpha v_i$ for $0 \leq i \leq d$.
\end{enumerate}
\end{theorem}

\proof Similar to the proof of Theorem \ref{thm:unique}. \hfill $\Box$\\

We now give various characterizations of a standard basis and dual standard basis.

\begin{theorem} \label{char1}
Let $\{w_i\}^{d}_{i=0}$ be a sequence of vectors in $W\!$, not all $0$.  Then $\{w_i\}^{d}_{i=0}$ is a standard basis for $W$ if and only if both (i) and (ii) below hold.
\begin{enumerate}
\item[\rm(i)]  $w_i \in E^{*}_{r+i}W$ $(0 \leq i \leq d)$.
\item[\rm(ii)] $ \sum^{d}_{i=0}w_i \in E_{t}W$\!.
\end{enumerate}
\end{theorem}

\proof Suppose that $\{w_i \}_{i=0}^d$ is a standard basis for $W$.  By Definition \ref{standardbasis}, there exists a nonzero $u \in E_tW$ such that $w_i = E_{r+i}^{*}u$ for $0 \leq i \leq d$.  Thus, (i) holds. Combining Lemma \ref{lemma:ds1}(ii) and the fact that $\sum_{j=0}^D E^*_{j} = I$, we have $u = \sum_{i=0}^dE^*_{r+i}u$. From this comment, we find that $\sum_{i=0}^d w_i = \sum_{i=0}^d E_{r+i}^{*}u = u \in E_tW$\!.  Hence, (ii) holds.  Conversely, suppose that $\{w_i \}_{i=0}^d$ satisfies (i) and (ii).  Let $u =  \sum^{d}_{i=0}w_i$.  By (ii) and the fact that not all of $\{w_i\}_{i=0}^d$ are $0$, $u$ is a nonzero vector in $E_tW$\!.  By (i), $E_{r+i}^{*}u = w_i$ for $0 \leq i \leq d$.  Therefore, $\{w_i\}^{d}_{i=0}$ is a standard basis for $W$.   \hfill $\Box$

\begin{theorem} \label{char1*}
Let $\{v_i\}^{d}_{i=0}$ be a sequence of vectors in $W\!$, not all $0$.  Then $\{v_i\}^{d}_{i=0}$ is a dual standard basis for $W$ if and only if both (i) and (ii) below hold.
\begin{enumerate}
\item[\rm(i)]  $v_i \in E_{t+i}W$ $(0 \leq i \leq d)$.
\item[\rm(ii)] $ \sum^{d}_{i=0}v_i \in E^*_{r}W$\!.
\end{enumerate}
\end{theorem}
\proof Similar to the proof of Theorem \ref{char1}. \hfill $\Box$

\begin{lemma} \label{char3}
Let $\{w_i \}^{d}_{i=0}$ be a basis for $W$\!.  With respect to this basis, let $B$ and $B^*$ denote the matrix representations of $A$ and $A^*$\!,  respectively.  Then $\{w_i \}^{d}_{i=0}$ is a standard basis for $W$ if and only if both (i) and (ii) below hold.
\begin{enumerate}
\item[\rm(i)] $B$ has constant row sum $\theta_t$.
\item[\rm(ii)] $B^* = \mbox{ diag}(\theta^{*}_r, \theta^{*}_{r+1}, \ldots, \theta^{*}_{r+d})$.
\end{enumerate}
\end{lemma}

\proof Let $w = \sum_{i=0}^d w_i$.  Note that $Aw = \sum_{i=0}^d\sum_{j=0}^d B_{ji}w_j$. Since $E_tW$ is the eigenspace of $A$ corresponding to $\theta_t$, by the previous statement, $B$ has constant row sum equal to $\theta_t$ if and only if $w \in E_tW$\!.  Observe also that $w_i \in E_{r+i}^{*}W$ if and only if $B^{*} = \mbox{ diag}(\theta^{*}_r, \theta^{*}_{r+1}, \ldots, \theta^{*}_{r+d})$.  The result follows from these comments and Theorem \ref{char1}.
\hfill $\Box$

\begin{lemma} \label{char3*}  Let $\{v_i \}^{d}_{i=0}$ be a basis for $W$\!.  With respect to this basis, let $B$ and $B^*$ be the matrix representations of $A$ and $A^*$\!, respectively.  Then $\{v_i \}^{d}_{i=0}$ is a dual standard basis for $W$ if and only if both (i) and (ii) below hold.
\begin{enumerate}
\item[\rm(i)] $B^{*}$ has constant row sum $\theta^{*}_r$.
\item[\rm(ii)] $B = \mbox{ diag}(\theta_t, \theta_{t+1}, \ldots, \theta_{t+d})$.
\end{enumerate}
\end{lemma}

\proof Similar to the proof of Lemma \ref{char3}.\hfill $\Box$

\begin{definition} \rm \label{maps}
Define the two maps $\flat: \mbox{ End}(W) \rightarrow \mbox{ Mat}_{d+1}(\mathbb{C})$ and $\sharp: \mbox{ End}(W) \rightarrow \mbox{ Mat}_{d+1}(\mathbb{C})$ as follows:  For every $Y \in \mbox{ End}(W), \ Y^{\flat}$ (resp.\! $Y^{\sharp}$) is the matrix representation of $Y$ with respect to a standard basis (resp.\! dual standard basis) for $W$\!.  Note that $Y^{\flat}$ (resp.\! $Y^{\sharp}$) is independent of the choice of standard basis (resp.\! dual standard basis) by Theorem \ref{thm:unique} (resp. \! Theorem \ref{thm:unique*}).  
\end{definition}

\begin{theorem}\label{sharpflat}
With reference to Definition \ref{maps}, the following hold. 
\begin{enumerate}
\item[\rm(i)] $A^{\flat}$ has constant row sum $\theta_t$.
\item[\rm(ii)] ${A^{*}}^{\flat} = \mbox{ diag}(\theta^{*}_r, \theta^{*}_{r+1}, \ldots, \theta^{*}_{r+d})$.
\item[\rm(iii)] ${A^{*}}^{\sharp}$ has constant row sum $\theta^{*}_r$.
\item[\rm(iv)] $A^{\sharp} = \mbox{ diag}(\theta_t, \theta_{t+1}, \ldots, \theta_{t+d})$.
\end{enumerate}
\end{theorem}

\proof Immediate from Lemmas \ref{char3} and \ref{char3*}. \hfill $\Box$

\section{The scalars $b_i(W), \ c_i(W)$} \label{seconbici}

Let $W$ be as in Assumption \ref{W} and let $\flat, \sharp$ be the maps in Definition \ref{maps}. In this section, we will take a close look at the entries of $A^{\flat}$ and $A^{*\sharp}$.  \\
\indent
By Lemmas \ref{eirep} and \ref{eirep*}, the matrices $A^{\flat}$ and $A^{*\sharp}$ are tridiagonal.  Moreover, by Lemma \ref{genbasis}, the $(i,i)$-entry of these matrices are $a_i(W)$ and $a^*_i(W)$, respectively. We now take a close look at the superdiagonal and subdiagonal entries of these matrices.  

\begin{definition} \rm\label{def:bi&ci}
Define
\begin{align*} 
b_i(W) \ & = \ (A^{\flat})_{i,i+1}, \qquad \qquad b^*_i(W) \ = \ (A^{*\sharp})_{i,i+1} \qquad  (0 \leq i \leq d-1),\\
c_i(W) \ & = \ (A^{\flat})_{i,i-1}, \qquad \qquad c^*_i(W) \  = \ (A^{*\sharp})_{i,i-1} \qquad (1 \leq i \leq d).
\end{align*}
Thus,
\begin{equation} A^{\flat} =  \left(
    \begin{array}{llcccrl} 
    a_0(W) &   b_0(W)   &          &        &                &      {\bf 0} \\
    c_1(W) &   a_1(W)   & b_1(W)   &        &                &              \\ 
           &   c_2(W)   & \cdot    &  \cdot &                &              \\ \label{matrixrepstandard} 
           &            & \cdot    &  \cdot &  \cdot         &              \\ 

           &            &          &  \cdot & a_{d-1}(W) & b_{d-1}(W)   \\
    {\bf 0}&            &          &        & c_d(W)         & a_d(W)  \\
    \end{array} \right),
\end{equation}    
\begin{equation} {A^{*}}^{\sharp} =  \left(
    \begin{array}{llcccrc}
    a_0^{*}(W) &   b_0^{*}(W)   &              &               &           &  {\bf 0} \\
    c_1^{*}(W) &   a_1^{*}(W)   & b_1^{*}(W)   &               &           &          \\
           &   c_2^{*}(W)       & \cdot        &  \cdot        &           &          \\
           &            & \cdot &  \cdot       &  \cdot        &          \\     \label{matrixrepstandard*}  
           &            &       &  \cdot       & a_{d-1}^{*}(W)& b_{d-1}^{*}(W) \\
    {\bf 0}&            &       &              & c_d^{*}(W)    & a_d^{*}(W)  \\
    \end{array}
\right). \end{equation}
For notational convenience, define $b_d(W)=0, \ c_0(W)= 0$ (resp.\! $b^*_d(W)=0,\ c^*_0(W) = 0$).  
Observe that by Lemmas \ref{eirep}(iii) and \ref{eirep*}(iii), $b_i(W), \ b^*_i(W) \ (0 \leq i \leq d-1), \ c_i(W), \ c^*_i(W) \ (1 \leq i \leq  d)$ are all nonzero. 
\end{definition}

\begin{definition} \rm
By the \textit{intersection numbers} (resp.\! \textit{dual intersection numbers}) of $W\!,$ we mean the $a_i(W)$, $b_i(W)$, $c_i(W)$ (resp.\! $ a_i^*(W), \ b^*_i(W), \ c^*_i(W)$).  
\end{definition}
     
 
\begin{lemma}\label{propbici}
The following hold.
\begin{enumerate}
\item[\rm(i)]  $b_{i-1}(W)c_i(W) = x_i(W) \quad (1 \leq i \leq d )$.
\item[\rm(ii)] $c_i(W) + a_i(W) + b_i(W) = \theta_{t} \quad (0 \leq i \leq d)$.
\item[\rm(iii)] $b^{*}_{i-1}(W)c^{*}_i(W) = x^{*}_i(W) \quad (1\leq i \leq d)$.
\item[\rm(iv)] $c^{*}_i(W) + a^{*}_i(W) + b^{*}_i(W) = \theta^{*}_{r} \quad (0 \leq i \leq d)$.
\item[\rm(v)] $b_i(W)\in \mathbb{R}, \ c_i(W) \in \mathbb{R} \quad (0 \leq i \leq d)$.
\item[\rm(vi)] $b_i^{*}(W)\in \mathbb{R}, \ c_i^{*}(W) \in \mathbb{R} \quad (0 \leq i \leq d)$.
\end{enumerate}
\end{lemma}

\proof (i)  Immediate from Lemma \ref{genbasis}(ii).\\
(ii) Immediate from Theorem \ref{sharpflat}(i).\\
(iii), (iv)  Similar to the proofs of (i), (ii).\\
(v) Recall that $a_0(W) \in \mathbb{R}$ by Lemma \ref{sumofai}(iii).  Since $\theta_t \in \mathbb{R}$ and $a_0(W) + b_0(W) = \theta_t$, we have $b_0(W) \in \mathbb{R}$.  By Lemma \ref{sumofai}(iii), (iv), we obtain $a_i(W) \in \mathbb{R}$ and $x_i(W) \in \mathbb{R}$ for $0 \leq i \leq d$.  Combining this with (i), (ii) and the fact that $b_0(W) \in \mathbb{R}$ and $b_i(W) \neq 0$ for $0 \leq i \leq d-1$, we find that $b_i(W)\in \mathbb{R}$ and $c_i(W) \in \mathbb{R}$ for $0 \leq i \leq d$.\\
(vi) Similar to the proof of (v). \hfill $\Box$

\begin{lemma} \label{bipi}
For $0 \leq i \leq d$, 
\begin{align}
b_0(W)b_1(W) \cdots b_{i-1}(W) \ &= \ p_i(\theta_t), \label{pibi}\\
b^{*}_0(W)b^{*}_1(W) \cdots b^{*}_{i-1}(W) \ &= \ p^{*}_i(\theta^{*}_r), \label{pibi*}
\end{align}
where $p_i = p_i^W, \ p^*_i = p^{*W}_i$ are from Definitions \ref{polypi}, \ref{polypi*}.
\end{lemma}

\proof We prove (\ref{pibi}) by induction on $i$. It can be verified that (\ref{pibi}) is true for $i=0, 1$.  Fix $2 \leq i \leq d$.  By (\ref{pi}),
\begin{equation}
p_i(\theta_t) =   (\theta_t - a_{i-1}(W))p_{i-1}(\theta_{t}) - x_{i-1}(W)p_{i-2}(\theta_t). \label{inductive2}
\end{equation}  
Eliminate $x_{i-1}(W)$ and $a_{i-1}(W)$ in (\ref{inductive2}) using Lemma \ref{propbici}(i), (ii).  Evaluate the result using the inductive hypothesis to obtain the desired result. Equation (\ref{pibi*}) is proved similarly.\hfill $\Box$
 
\begin{theorem} \label{formbici}
The following (i)--(iv) hold.
\begin{enumerate}
\item[{\rm(i)}] $b_{i}(W) = {\displaystyle\frac{p_{i+1}(\theta_{t})}{p_i(\theta_{t})}\quad  (0 \leq i \leq d-1)}$. \\
\item[{\rm(ii)}] $c_{i}(W) = {\displaystyle\frac{x_{i}(W)p_{i-1}(\theta_{t})}{p_{i}(\theta_{t})} \quad (1 \leq i \leq d)}$.\\
\item[{\rm(iii)}] $b^{*}_{i}(W) = {\displaystyle\frac{p^{*}_{i+1}(\theta^*_{r})}{p^{*}_i(\theta^{*}_{r})}  \quad (0 \leq i \leq d-1)}$. \\
\item[{\rm(iv)}] $c^{*}_{i}(W) = {\displaystyle\frac{x^{*}_{i}(W)p^{*}_{i-1}(\theta^*_{r})}{p^{*}_{i}(\theta^{*}_{r})} \quad (1 \leq i \leq d)}$.
\end{enumerate}
In the above lines, $p_j = p_j^W, \ p^*_j = p^{*W}_j$ are from Definitions \ref{polypi}, \ref{polypi*}.
\end{theorem}

\proof (i)  Immediate from Lemma \ref{bipi}.\\
(ii) Immediate from (i) and Lemma \ref{propbici}(i).\\
(iii), (iv) Similar to the proofs of (i), (ii).\hfill $\Box$

\begin{lemma}{\rm\cite[Theorem 4.1(vi)]{subII}}\label{fortrivial}
Let $W$ be the trivial $T$-module.  For $0 \leq i \leq D$, let $a_i, \ b_i, \ c_i$ (resp. $a^*_i, \ b^*_i, \ c^*_i$) be the intersection (resp. dual intersection) numbers of $\Gamma$.  Then 
\begin{enumerate}
\item[\rm(i)] $a_i(W) = a_i, \ b_i(W) = b_i, \ c_i(W) = c_i$,
\item[\rm(ii)] $a^*_i(W) = a^*_i, \ b^*_i(W) = b^*_i, \ c^*_i(W) = c^*_i$\!.
\end{enumerate}
\end{lemma}

\indent We finish this section with a few comments.

\begin{lemma} \label{isoclass}
Let $W, \ W'$ be thin irreducible $T\!$-modules.  The following are equivalent.
\begin{enumerate}
\item[\rm(i)] $W$ and $W'$ are isomorphic $T\!$-modules.
\item[\rm(ii)] $W$ and $W'$ have the same endpoint, dual endpoint, diameter and intersection numbers.
\item[\rm(iii)] $W$ and $W'$ have the same endpoint, dual endpoint, diameter and dual intersection numbers.
\end{enumerate}   
\end{lemma}

\proof
(i) $\Rightarrow$ (ii) Suppose that $W$ and $W'$ are isomorphic $T\!$-modules.  Let $\phi: W \rightarrow W'$ be an isomorphism of $T\!$-modules.  Thus, $\phi(E_iW) = E_iW'$.  Hence $E_iW \neq 0 $ if and only if $E_iW' \neq 0$.  Similarly, $E^{*}_iW \neq 0 $ if and only if $E^{*}_iW' \neq 0$.  Therefore, $W$ and $W'$ have the same endpoint, dual endpoint and diameter.  Since $W$ and $W'$ are isomorphic, the matrices representing the action of $A$ on $W$ and $W'$ are the same.  Hence they have the same intersection numbers.    \\  
(ii) $\Leftarrow$ (i) Suppose that $W$ and $W'$ have the same endpoint $r$, dual endpoint $t$ and diameter $d$.  Suppose also that they have the same intersection numbers.  For $0 \leq i \leq d$, let $w_i = E^*_{r+i}u$ and $ w_i' = E^*_{r+i}u'$, where $u$ and $u'$ are nonzero vectors in $E_tW$ and $E_tW'\!,$ respectively.  Since $W$ and $W'$ both have dimension $d+1$, there exists a vector space isomorphism $\phi: W \rightarrow W'$ such that $\phi(w_i) = w_i'$.   Since $w_i \in E^*_{r+i}W$, it can be easily verified that $(\phi A^* - A^*\phi)w_i = 0$ for $0 \leq i \leq D$.  By (\ref{matrixrepstandard}) and the fact that $W$ and $W'$ have the same intersection numbers, $(\phi A - A\phi)w_i = 0$ for $0 \leq i \leq d$.  From these comments, $(\phi A - A\phi)W = 0$ and $(\phi A^* - A^*\phi)W=0$.  Since $T$ is generated by $A, A^*$, we find that $\phi$ is a $T$-module isomorphism.  Therefore, $W$ and $W'$ are isomorphic $T$-modules.\\
(i) $\Leftrightarrow$ (iii) Similar to the proof of (i) $\Leftrightarrow$ (ii). \hfill $\Box$

\section {The scalar $k_i(W)$}

Let $W$ be as in Assumption \ref{W}.  In this section, we will look at a sequence of scalars closely related with the $m_i(W)$.

\begin{definition} \rm \label{ki}
For $0 \leq i \leq d$, define
\begin{align*}
k_i(W) \ &= \ m^{*}_i(W)\nu(W),\\
k^{*}_i(W) \ &= \  m_i(W)\nu(W),
\end{align*}
where $m_i(W), \ m^*_i(W), \ \nu(W)$ are from Definitions \ref{defmi} and \ref{defnu}.

\end{definition}

\begin{lemma} \label{propki}
The following (i)--(iii) hold.
   \begin{enumerate}
   \item[\rm(i)] $k_0(W) = 1, \ k_0^{*}(W) = 1$.
    \item[\rm(ii)] $\sum_{i=0}^{d} k_i(W) = \nu(W)$.
    \item[\rm(iii)] $ \sum_{i=0}^{d} k_i^{*}(W) = \nu(W)$.
    \item[\rm(iv)] $k_i(W) >0, \ k^*_i(W) >0  \quad (0 \leq i \leq d)$.
   \end{enumerate}
In the above lines, $\nu(W)$ is from Definition \ref{defnu}.   
\end{lemma}
\proof (i) Immediate from Definition \ref{ki}.\\
(ii) Immediate from Lemma \ref{miprop2}(i) and Definition \ref{ki}.\\
(iii) Similar to the proof of (ii).\\
(iv) Immediate from Lemma \ref{miprop2}(iii), (iv) and Definitions \ref{defnu}, \ref{ki}. \hfill $\Box$\\

We now relate $k_i(W)$ (resp.\! $k^*_i(W)$) and the intersection (resp.\! dual intersection) numbers of $W$. 

\begin{lemma} \label{lemma:kirec}
For $0 \leq i \leq d$,
\begin{align}
k_i(W)c_i(W) \ & =  \ k_{i-1}(W)b_{i-1}(W),\label{recki}\\
k^*_i(W)c^*_i(W) \ & =  \ k^*_{i-1}(W)b^*_{i-1}(W), \label{recki*}
\end{align}
where $b_j(W),\  b^*_j(W), \ c_j(W), \ c^*_j(W)$ are from Definition \ref{def:bi&ci} and $b_{-1}(W) =0, \ b^*_{-1}(W) =0 $.
\end{lemma}

\proof
We proceed by induction on $i$.
Since $c_0(W) = 0$, equation (\ref{recki}) holds for $i=0$. Assume $1 \leq i \leq d$.  By Definition \ref{def:bi&ci}, on $W$
\begin{equation}
AE^*_{r+i}E_t = b_{i-1}(W)E^*_{r+i-1}E_t + a_{i}(W)E^*_{r+i}E_t + c_{i+1}(W)E^*_{r+i+1}E_{t},\label{eqn1}
\end{equation}
where $c_{d+1}(W)= 0$.
Take the trace of both sides of (\ref{eqn1}). Evaluate this using Definition \ref{defmi} and the fact that $E_tA = \theta_tE_t$.  Multiplying $\nu(W)$ on both sides of the resulting equation and using Definition \ref{ki} we obtain
\begin{equation}
\theta_t k_{i}(W) = b_{i-1}(W)k_{i-1}(W) + a_{i}(W)k_{i}(W) + c_{i+1}(W)k_{i+1}(W). \label{eqki}
\end{equation}
Solving for $c_{i+1}(W)k_{i+1}(W)$ in (\ref{eqki}) using the inductive hypothesis and Lemma \ref{propbici}(ii), we find that (\ref{recki}) holds for $i+1$. The proof of (\ref{recki*}) is similar. \hfill $\Box$

\begin{theorem}\label{kiintermsofbi}
For $0 \leq i \leq d$,
\begin{align}
k_i(W) \ &= \ \frac{b_0(W)b_1(W)\cdots b_{i-1}(W)}{c_1(W)c_2(W)\cdots c_{i}(W)}, \label{kiform}\\\
k^{*}_i(W) \ &= \ \frac{b^{*}_0(W)b^*_1(W)\cdots b^{*}_{i-1}(W)}{c^{*}_1(W)c^*_2(W)\cdots c^{*}_{i}(W)}, \label{kidualform}
\end{align}
where $b_j(W),\  b^*_j(W), \ c_j(W), \ c^*_j(W)$ are from Definition \ref{def:bi&ci}.
\end{theorem}

\proof Solve for $k_i(W)$ and $k^*_i(W)$ in Lemma \ref{lemma:kirec} recursively to obtain the desired result.  \hfill $\Box$   

\begin{corollary}
Let $W$ be the trivial $T$-module.  Then for $0 \leq i \leq D$,
$$k_i(W) = k_i, \qquad \qquad k^*_i(W) = m_i,$$
where $k_i$ is the $i$th valency of $\Gamma$ and $m_i$ is the multiplicity of $\Gamma$ associated with $E_i$.
\end{corollary}

\proof
Immediate from (\ref{valencyform}), (\ref{multform}), Lemma \ref{fortrivial} and Theorem \ref{kiintermsofbi}. \hfill $\Box$



\section{The polynomials $u_i$ and $v_i$}

Let $W$ be as in Assumption \ref{W}.  In this section, we will look at two normalizations of the polynomials $p_i$ and $p_i^*$ in Definitions \ref{polypi}, \ref{polypi*}.  

\begin{definition} \rm\label{defdefvi}
Define $v_i = v_i^{W} \mbox{ and } v^{*}_i = v^{*W}_i $ in $\mathbb{C}[\lambda]$ by
\begin{align}
v_i  \ &= \  \frac{p_i}{c_1(W)c_2(W) \cdots c_i(W)} \qquad(0 \leq i \leq d),  \label{defvi}\\
v^{*}_i  \ &= \  \frac{p^{*}_i}{c^{*}_1(W)c^{*}_2(W) \cdots c^{*}_i(W)} \qquad(0 \leq i \leq d),
\end{align}
where $p_i = p_i^W, \ p^*_i = p^{*W}_i$ are from Definitions \ref{polypi}, \ref{polypi*} and $ c_j(W), \ c^*_j(W)$ are from Definition \ref{def:bi&ci}.
For notational convenience, define $v_{-1} = 0, \ v^{*}_{-1}=0 $.
\end{definition}

\begin{lemma}\label{viki}
For $0 \leq i \leq d$,   
$$v_i(\theta_t) = k_i(W), \qquad \qquad v^{*}_i(\theta^{*}_r) = k^{*}_i(W),$$
where 
$v_i = v_i^W$, $v^*_i=v^{*W}_i$ are from Definition \ref{defdefvi}. 
\end{lemma}

\proof Immediate from Lemma \ref{bipi}, Theorem \ref{kiintermsofbi} and Definition \ref{defdefvi}.   \hfill $\Box$

\begin{lemma}  With reference to Definition \ref{defdefvi}, for $0 \leq i \leq d-1$,  
\begin{align}
\lambda v_i \ &= \  b_{i-1}(W)v_{i-1}+ a_i(W)v_i + c_{i+1}(W)v_{i+1}, \label{eq:vi}\\
\lambda v^*_i \ &= \  b^*_{i-1}(W)v^*_{i-1}+ a^*_i(W)v^*_i + c^*_{i+1}(W)v^*_{i+1},\label{eq:vi*}
\end{align}
where $b_{-1}(W) = 0, \ b^*_{-1}(W) = 0$.
Moreover,
\begin{align*}
\lambda v_d - a_d(W)v_d - b_{d-1}(W)v_{d-1} \ &= \ c^{-1}p_{d+1},\\
\lambda v^{*}_d - a^{*}_d(W)v^{*}_d - b^{*}_{d-1}(W)v^{*}_{d-1} \ &= \ c^{*-1}p^{*}_{d+1},
\end{align*}
where
\begin{align*}
c \ &= \ c_1(W)c_2(W)\cdots c_d(W),\\
c^* \ &= \ c^{*}_1(W)c^{*}_2(W)\cdots c^{*}_d(W).
\end{align*}
\end{lemma}

\proof To obtain (\ref{eq:vi}), divide both sides of (\ref{pi}) by $c_1(W)c_2(W)\cdots c_i(W)$ and eliminate $x_i(W)$ using Lemma \ref{propbici}(i). The proof of (\ref{eq:vi*}) is similar.\hfill $\Box$
\begin{theorem} \label{vpush}

With reference to Definition \ref{defdefvi}, for $0 \leq i \leq d$,
\begin{equation}
v_i(A)E^{*}_{r}u =  E^{*}_{r+i}u, \qquad \qquad v^{*}_i(A^{*})E_{t}v = E_{t+i}v, \label{viEi}
\end{equation}
where $u$ and $v$ are nonzero vectors in  $E_tW$ and $E^*_rW$, respectively.
\end{theorem}

\proof For $0 \leq i \leq d$, let $w_i = E_{r+i}^{*}u \mbox{ and }  w_i' = v_i(A)E_r^{*}u $. By (\ref{matrixrepstandard}),
\begin{equation}
Aw_i = b_{i-1}(W)w_{i-1} + a_i(W)w_i + c_{i+1}(W)w_{i+1} \quad (0 \leq i \leq d-1), \label{wirec}
\end{equation}
where $b_{-1}(W)=0$.
Using (\ref{eq:vi}), we obtain
\begin{equation}
Aw_i' = b_{i-1}(W)w_{i-1}' + a_i(W)w_i' + c_{i+1}(W)w_{i+1}' \quad (0 \leq i \leq d-1). \label{wi'rec}
\end{equation}
Using the fact that $w_0 = w_0'$ and comparing (\ref{wirec}) and (\ref{wi'rec}), we obtain the equation on the left of (\ref{viEi}). The equation on the right of (\ref{viEi}) can be similarly obtained. \hfill $\Box$

\begin{definition} \rm \label{defdefui}
For $0 \leq i \leq d$, define $u_i= u_i^{W}$ and $u^{*}_i= u^{*W}_i$ in $\mathbb{C}[\lambda]$ as follows:
\begin{align}
u_i  \ &= \  \frac{p_i}{p_i(\theta_t)}, \label{defui}\\
u^{*}_i  \ &= \  \frac{p^{*}_i}{p^*_i(\theta^*_r)},
\end{align}
where  $p_i = p_i^W, \ p^*_i = p_i^{*W}$ are from Definitions \ref{polypi}, \ref{polypi*}.
For notational convenience, define $u_{-1}=0, \ u^*_{-1}=0$.
\end{definition}

\begin{lemma} \label{viuiki}
With reference to Definition \ref{defdefvi}, for $0 \leq i \leq d$,
$$v_i = k_i(W)u_i, \qquad \qquad v^{*}_i = k^{*}_i(W)u^{*}_i,$$
where $u_i = u_i^{W}, \ u^*_i = u_i^{*W}$ are from Definition \ref{defdefui} and $k_i(W), \ k^*_i(W)$ are from Definition \ref{ki}. 
\end{lemma}

\proof Immediate from Lemma \ref{bipi}, Theorem \ref{kiintermsofbi} and Definitions \ref{defdefvi}, \ref{defdefui}. \hfill $\Box$

\begin{lemma} \label{recurrenceui}
With reference to Definition \ref{defdefui}, for $0 \leq i \leq d-1$,
\begin{align}
\lambda u_i \ &= \ c_i(W)u_{i-1} + a_i(W)u_i + b_{i}(W)u_{i+1},  \label{recforui}\\
\lambda u^{*}_i \ &= \ c^{*}_i(W)u^{*}_{i-1} + a^{*}_i(W)u^{*}_i + b^{*}_{i}(W)u^{*}_{i+1}. \label{recforui*}
\end{align}
Moreover,
\begin{align*}
\lambda u_d - c_d(W)u_{d-1} - a_d(W)u_{d} \ &= \ p_{d+1}/p_d(\theta_t),\\
\lambda u^{*}_d - c^{*}_d(W)u^{*}_{d-1} - a^{*}_d(W)u^{*}_{d} \ &= \ p^{*}_{d+1}/p^{*}_d(\theta^{*}_r).
\end{align*}
\end{lemma}
\proof 
To obtain (\ref{recforui}), divide both sides of (\ref{pi}) by $p_i(\theta_t)$ and eliminate $x_i(W)$ using Lemma \ref{propbici}(i).  Evaluate the result using Lemma \ref{bipi}.  The proof of (\ref{recforui*}) is similar.  \hfill $\Box$

\begin{theorem}\label{3term}
With reference to Definition \ref{defdefui}, for $0 \leq i,j \leq d$,
\begin{align*}
\theta_{t+j}u_i(\theta_{t+j}) \ &= \ c_{i}(W)u_{i-1}(\theta_{t+j}) + a_{i}(W)u_{i}(\theta_{t+j}) + b_{i}(W)u_{i+1}(\theta_{t+j}),\\
\theta^{*}_{r+j}u^{*}_i(\theta^{*}_{r+j}) \ &= \  c^{*}_{i}(W)u^{*}_{i-1}(\theta^{*}_{r+j})  + a^{*}_{i}(W)u^{*}_{i}(\theta^{*}_{r+j}) + b^{*}_{i}(W)u^{*}_{i+1}(\theta^{*}_{r+j}),
\end{align*}
where $u_{d+1}=0, \ u^*_{d+1}=0.$
\end{theorem}

\proof
Immediate from (\ref{recforui}) and (\ref{recforui*}) with $\lambda = \theta_{t+j} \mbox{ and } \lambda = \theta^{*}_{r+j}$. \hfill $\Box$

\section  {Some inner products and the Askey-Wilson duality}
Let $W$ be as in Assumption \ref{W}.   In this section, we will look at all inner products involving the elements of a standard basis and a dual standard basis for $W$.  Using these inner products, we will show that all the polynomials associated with $W$ satisfy relations known as the Askey-Wilson duality.

Throughout the entire section, $u$ and $v$ are nonzero vectors in $E_tW\!$  and $E^*_rW$\!, respectively. 
Recall that by Definition \ref{standardbasis}, $\{E^*_{r+i}u \}_{i=0}^d$ (resp.\! $\{E_{t+i}v \}_{i=0}^d$) is a standard basis (resp.\! dual standard basis) for $W$.  By (\ref{eq:decomposeV}) and (\ref{eq:decomposeV2}), each of these bases is orthogonal.  We now compute some square norms. 

\begin{theorem} \label{innersame}
For $0 \leq i \leq d$,
\begin{eqnarray}
\|E^{*}_{r+i}u\|^2& = & \|u\|^{2} k_i(W)  /\nu(W), \label{innersame1}\\
\|E_{t+i}v\|^2 & = &  \|v\|^{2}k^{*}_i(W)/\nu(W), \label{innersame2}
\end{eqnarray}
where $\nu(W)$ is from Definition \ref{defnu} and $k_i(W), \ k^*_i(W)$ are from Definition \ref{ki}.
\end{theorem}

\proof Note that
\begin{eqnarray*}
\|E^{*}_{r+i}u\|^2 & = & \langle E^{*}_{r+i}u,E^{*}_{r+i}u \rangle \\
 & = & \langle u, E^{*2}_{r+i}u\rangle\\
& = & \langle u, E^{*}_{r+i}u\rangle\\
& = & \langle u, v_i(A)E_r^{*}u\rangle \qquad \quad \ \mbox{ by Lemma \ref{vpush}}\\
& = & \langle v_i(A)u, E_r^{*}u\rangle\\
& = &\langle v_i(\theta_t)u, E_r^{*}u\rangle\\
& = & k_i(W)\langle u, E_r^{*}u\rangle \qquad \quad\mbox{ by Lemma \ref{viki}}.
\end{eqnarray*}
Since $u \in E_tW\!$, $u = E_tu$. Using this we find that $\langle u, E_r^{*}u \rangle = \langle E_tu, E^*_rE_tu \rangle = \langle u, E_tE^*_rE_tu \rangle$.  Evaluating $E_tE^*_rE_t$ using Lemma \ref{thm:nu}(i) we find that $\langle u, E_r^{*}u \rangle = \|u \|^2/\nu(W)$.  Thus, we obtain (\ref{innersame1}). Equation (\ref{innersame2}) is proved similarly. \hfill $\Box$\\

\indent
Our next goal is to compute the inner product between the elements of $\{E^*_{r+i}u \}_{i=0}^d$ and $\{E_{t+i}v \}_{i=0}^d$.  We need the following lemma. 

\begin{lemma} \label{er*et}
The following hold.
\begin{enumerate}
\item[\rm(i)] $\langle E^{*}_ru, E_tv \rangle = \langle u, v \rangle/\nu(W)$.
\item[\rm(ii)]$E^{*}_r u = \frac{\langle u, v \rangle}{\|v \|^{2}}v$.
\item[\rm(iii)] $E_tv = \frac{\langle v, u \rangle}{\|u \|^{2}}u$. 
\item[\rm(iv)] $\langle u, v \rangle \neq 0$.
\item[\rm(v)] $\nu(W) |\langle u,v \rangle|^2 = \| u\|^2 \|v \|^2$.
\end{enumerate}
In the above lines, $\nu(W)$ is from Definition \ref{defnu}.
\end{lemma}

\proof (i) Since $v \in E^*_rW\!,$ $v=E^*_rv$.  Using this we find that $\langle E^{*}_ru, E_tv \rangle = \langle E^*_ru, E_tE_r^{*}v \rangle = \langle u, E^*_rE_tE_r^{*}v \rangle $.  Evaluate $E^*_rE_tE^*_r$ using Lemma \ref{thm:nu}(ii) to obtain the desired result.\\
(ii) Since $v$ spans $E_r^{*}W$\!, $E_r^{*}u = \alpha v$ for some $\alpha \in \mathbb{C}$.  Thus $\langle E_r^{*}u, v \rangle = \alpha \|v \|^2$.  Since $\langle E_r^{*}u, v \rangle = \langle u, E^*_rv \rangle = \langle u, v \rangle$, we find that $\alpha = \frac{\langle u, v\rangle}{\| v\|^2} $.\\
(iii) Similar to the proof of (ii).\\
(iv) Observe that $E^*_ru \neq 0$ since it is an element of a standard basis.  It follows from this and (ii) that $\langle u, v \rangle \neq 0$.\\  
(v) Eliminate $E^{*}_ru$ and $E_tv$ in (i) using (ii) and (iii).\hfill $\Box$

\begin{theorem}\label{innerdif}
For $0 \leq i,j \leq d$,
\begin{align}
\langle E^{*}_{r+i}u, E_{t+j}v \rangle  \ &= \  u_i(\theta_{t+j})  k_{i}(W)k^{*}_{j}(W)\langle u,v\rangle/\nu(W), \label{old}\\
\langle E^{*}_{r+i}u, E_{t+j}v \rangle\ &= \  u_j^{*}(\theta_{r+i}^{*}) k_{i}(W)k^{*}_{j}(W) \langle u,v\rangle /\nu(W), \label{new}
\end{align}
where $\nu(W), \ k_i(W), \ k^*_j(W)$ are from Definitions \ref{defnu}, \ref{ki} and \ $u_i=u_i^W\!, \ u^*_j = u_j^{*W}$ are from Definition \ref{defdefui}.
\end{theorem}

\proof Note that
\begin{align*}
\langle E^{*}_{r+i}u, E_{t+j}v \rangle \ & = \ \langle v_{i}(A)E_{r}^{*}u, E_{t+j}v \rangle & & \mbox{ by Theorem \ref{vpush} }\\
\ & = \ \langle E_{r}^{*}u,  v_{i}(A)E_{t+j}v \rangle\\
\ & = \ v_{i}(\theta_{t+j})\langle E_{r}^{*}u,  E_{t+j}v \rangle \\
\ & = \ v_{i}(\theta_{t+j}) \langle E_{r}^{*}u,  v_{j}^{*}(A^{*})E_{t}v \rangle & & \mbox{ by Theorem \ref{vpush}}\\
\ & = \ v_{i}(\theta_{t+j}) \langle v_{j}^{*}(A^{*})E_{r}^{*}u,  E_{t}v \rangle\\
\ & = \ v_{i}(\theta_{t+j})v_{j}^{*}(\theta_r^{*}) \langle E_{r}^{*}u,  E_{t}v \rangle\\
\ & = \ v_{i}(\theta_{t+j})v_{j}^{*}(\theta_r^{*}) \langle u, v \rangle/\nu(W) & & \mbox{ by Lemma \ref{er*et}(i)}.
\end{align*}
The result then follows from Lemmas \ref{viki} and \ref{viuiki}. Equation (53) is proved similarly. \hfill $\Box$


\begin{theorem}\label{askey}
For $0 \leq i, j \leq d$,
\begin{equation}
u_i(\theta_{t+j}) = {u^{*}_j}(\theta^{*}_{r+i}), \label{AWforui}
\end{equation}
where $u_i= u^W_i$ and $u^*_j = u_j^{*W}$ are from Definition \ref{defdefui}.
\end{theorem}

\proof
Compare (\ref{old}) with (\ref{new}).\hfill $\Box$

\begin{theorem}\label{thm:askey2}
For $0 \leq i, j \leq d$,
\begin{align}
\frac{p_i(\theta_{t+j})}{p_i(\theta_t)} \ & = \ \frac{p^{*}_j(\theta^{*}_{r+i})}{{p^{*}_j}(\theta^{*}_r)},  \label{AWDpi}\\
\frac{v_i(\theta_{t+j})}{k_i(W)} \ &= \ \frac{{v^{*}_j}(\theta^{*}_{r+i})}{k^{*}_j(W)},  \label{AWDvi}
\end{align}
where $p_i=p_i^W\!, \ p^*_i=p_i^{*W}\!, \ v_i =v_i^W\!, \ v^*_i =v_i^{*W}$ are from Definitions \ref{polypi}, \ref{polypi*} and \ref{defdefvi}.
\end{theorem}

\proof Immediate from Definition \ref{defdefui} and Theorems \ref{viuiki}, \ref{askey}. \hfill $\Box$\\

Equations (\ref{AWforui}), (\ref{AWDpi}) and (\ref{AWDvi}) are known as the \textit{Askey-Wilson duality}.  Combining Theorem \ref{3term} and Theorem \ref{askey}, we obtain the following result.

\begin{theorem}
For $0 \leq i, j \leq d$,
\begin{align}
\theta_{t+j}u^{*}_j(\theta^{*}_{r+i}) \ & = \ b_{i}(W)u^{*}_{j}(\theta^{*}_{r+i+1}) + a_{i}(W)u^{*}_{j}(\theta^{*}_{r+i}) + c_{i}(W)u^{*}_{j}(\theta^{*}_{r+i-1}), \label{difeqn1}\\
\theta^*_{r+j}u_j(\theta_{t+i})  \ & = \ b^{*}_{i}(W)u_{j}(\theta_{t+i+1}) + a^{*}_{i}(W)u_{j}(\theta_{t+i}) + c^{*}_{i}(W)u_{j}(\theta_{t+i-1}),\label{difeqn2}
\end{align}
where $u_j = u_j^W$ and $u^*_j = u_j^{*W}$ are from Definition \ref{defdefui}.
\end{theorem}

\section{The orthogonality relations}

Let $W$ be as in Assumption \ref{W}.  In this section, we display the transition matrix relating a standard basis and a dual standard basis.  Using this and the results of the previous section, we display the orthogonality relations satisfied by the polynomials we have seen in this paper. 

\begin{theorem} 
\label{transmatrix}
Let $u$ and $v$ be nonzero vectors in $E_tW$ and $E^*_rW\!$, respectively. For $0 \leq i \leq d$,  
\begin{align}
E^{*}_{r+i}u \ &= \ {\displaystyle \frac{\langle u,  v \rangle}{\| v\|^2} \sum_{j=0}^{d} v_i(\theta_{t+j})E_{t+j}v}, \label{trans1}\\
E_{t+i}v  \ &= \ {\displaystyle \frac{\langle v,  u \rangle}{\| u\|^2} \sum_{j=0}^{d} v^{*}_i(\theta^{*}_{r+j})E^{*}_{r+j}u},\label{trans2}
\end{align}
where $v_i = v_i^W\!, \ v^*_i = v_i^{*W}$ are from Definition \ref{defdefvi}.
\end{theorem}

\proof
Combining Lemma \ref{lemma:ds2}(ii) and the fact that $\sum_{j=0}^D E_{j} = I$, we find that $v = \sum_{j=0}^dE_{t+j}v$.  By Theorem \ref{vpush} and Lemma \ref{er*et}(ii), $E^*_{r+i}u = \frac{\langle u,v\rangle}{\|v \|^2}v_i(A)v$. Therefore, 
\begin{align*}
E^*_{r+i}u \ &= \ \frac{\langle u,v\rangle}{\|v \|^2}v_i(A)v \\
\ &= \ \frac{\langle u,v\rangle}{\|v \|^2}v_i(A)\sum_{j=0}^dE_{t+j}v\\
\ &= \ \frac{\langle u,  v \rangle}{\| v\|^2} \sum_{j=0}^{d} v_i(\theta_{t+j})E_{t+j}v.
\end{align*}
Hence, (\ref{trans1}) holds.  Equation (\ref{trans2}) is proved similarly. \hfill $\Box$   

\begin{theorem}\label{orthovi}
For $0 \leq i,j \leq d$, 
\begin{eqnarray}
{\displaystyle \sum_{h=0}^{d} v_i(\theta_{t+h})v_j(\theta_{t+h}) k^{*}_h(W)} & = & \delta_{ij}\nu(W) k_i(W), \label{ortho1}  \\
{\displaystyle \sum_{h=0}^{d} v_h(\theta_{t+i})v_h(\theta_{t+j})(k_h(W))^{-1}} & = & \delta_{ij}\nu(W) (k^{*}_i(W))^{-1} \label{ortho2}
\end{eqnarray}
and
\begin{eqnarray}
{\displaystyle \sum_{h=0}^{d} v^{*}_i(\theta^{*}_{r+h})v^*_j(\theta^{*}_{r+h}) k_h(W)} & = & \delta_{ij}\nu(W) k^{*}_i(W),\label{ortho9}\\
{\displaystyle \sum_{h=0}^{d} v^{*}_h(\theta^{*}_{r+i})v^{*}_h(\theta^{*}_{r+j})(k^{*}_h(W))^{-1}} & = & \delta_{ij}\nu(W) (k_i(W))^{-1}, \label{ortho10}
\end{eqnarray}
where $\nu(W), \ k_h(W), \ k_h^{*}(W)$ are from Definitions \ref{defnu}, \ref{ki} and $v_h= v_h^W\!, \ v^*_h= v^{*W}_h$ are from Definition \ref{defdefvi}.
\end{theorem}

\proof Concerning (\ref{ortho1}), let $u$ be a nonzero vector in $E_tW\!.$  We compute $\langle E_{r+i}^{*}u, E_{r+j}^{*}u \rangle$ in two ways.  First, by (\ref{eq:decomposeV2}) and (\ref{innersame1}), $\langle E_{r+i}^{*}u, E_{r+j}^{*}u \rangle = \delta_{ij}\|u\|^2k_i(W)/\nu(W)$.  Secondly, we compute $\langle E_{r+i}^{*}u, E_{r+j}^{*}u \rangle$ by evaluating each of $E_{r+i}^{*}u$ and $E_{r+j}^{*}u$  using (\ref{trans1}).  Simplify the result using (\ref{innersame2}) and Lemma \ref{er*et}(v).  We find that $\langle E_{r+i}^{*}u, E_{r+j}^{*}u \rangle$ is equal to $\|u\|^2/(\nu(W))^{2}$ times the left side of (\ref{ortho1}). Equation (\ref{ortho1}) follows from these comments. Similarly, we obtain (\ref{ortho9}).  To obtain (\ref{ortho2}), evaluate (\ref{ortho9}) using (\ref{AWDvi}). To obtain (\ref{ortho10}), evaluate (\ref{ortho1}) using (\ref{AWDvi}).\hfill $\Box$

\begin{theorem}\label{orthoui}
For $0 \leq i, j \leq d$,
\begin{eqnarray}
{\displaystyle \sum_{h=0}^{d} u_i(\theta_{t+h})u_j(\theta_{t+h}) k^{*}_h(W)} & = & \delta_{ij}\nu(W) (k_i(W))^{-1}, \label{ortho3}\\
{\displaystyle \sum_{h=0}^{d} u_h(\theta_{t+i})u_h(\theta_{t+j})k_{h}(W)} & = & \delta_{ij}\nu(W) (k^{*}_i(W))^{-1}, \label{ortho4}
\end{eqnarray}
and
\begin{eqnarray}
{\displaystyle \sum_{h=0}^{d} u^{*}_i(\theta^{*}_{r+h})u^{*}_j(\theta^{*}_{r+h}) k_h(W)} & = & \delta_{ij}\nu(W) (k_i^{*}(W))^{-1},\label{ortho7}\\
{\displaystyle \sum_{h=0}^{d} u^{*}_h(\theta^{*}_{r+i})u^{*}_h(\theta^{*}_{r+j})k^{*}_{h}(W)} & = & \delta_{ij}\nu(W) (k_i(W))^{-1},  \label{ortho8}
\end{eqnarray}
where $\nu(W), \ k_h(W), \ k_h^{*}(W)$ are from Definitions \ref{defnu}, \ref{ki} and $u_h= u_h^W\!, \ u^*_h= u^{*W}_h$ are from Definition \ref{defdefui}.
\end{theorem}

\proof 
Evaluate each of (\ref{ortho1})-(\ref{ortho10}) using Lemma \ref{viuiki}.\hfill $\Box$ 

\begin{theorem} \label{orthopi}
For $0 \leq i, j \leq d$,  
\begin{eqnarray}
{\displaystyle \sum_{h=0}^{d} p_i(\theta_{t+h})p_j(\theta_{t+h}) k^*_h(W)} & = & \delta_{ij}\nu(W)x_1(W)x_2(W)\cdots x_i(W),  \label{ortho5}\\
{\displaystyle \sum_{h=0}^{d} \frac{p_h(\theta_{t+i})p_h(\theta_{t+j})}{x_1(W)x_2(W) \cdots x_h(W)}} & = & \delta_{ij}\nu(W)(k^*_i(W))^{-1},\label{ortho6}
\end{eqnarray}
and
\begin{eqnarray}
{\displaystyle \sum_{h=0}^{d} p^{*}_i(\theta^{*}_{r+h})p^{*}_j(\theta^{*}_{r+h}) k_h(W)} & = & \delta_{ij}\nu(W) x^{*}_1(W)x^{*}_2(W)\cdots x^{*}_i(W),  \label{ortho11}\\
{\displaystyle \sum_{h=0}^{d} \frac{p^{*}_h(\theta^{*}_{r+i})p^{*}_h(\theta^{*}_{r+j})}{x^{*}_1(W)x^{*}_2(W) \cdots x^{*}_h(W)}} & = & \delta_{ij}\nu(W)(k_i(W))^{-1}, \label{ortho12}
\end{eqnarray}
where $ x_h(W),\ \nu(W), \ k_h(W), \ k_h^{*}(W)$ are from Definitions \ref{defxiai}, \ref{defnu}, \ref{ki} and $p_h = p_h^W\!, \ p^*_h = p_h^{*W} $ are from Definitions \ref{polypi}, \ref{polypi*}.
\end{theorem}

\proof Evaluate each of (\ref{ortho1})-(\ref{ortho10}) using Definition \ref{defdefvi}.  Simplify the result using Lemma \ref{propbici}(i), (iii). \hfill $\Box$\\

\indent
We now present Theorem \ref{orthovi} in matrix form.

\begin{definition} \rm\label{pmatrix}
Define matrices $P = P(W)$ and $P^{*}= {P^{*}}(W)$ in $\mbox{\rm{Mat}}_{d+1}(\mathbb{C})$ as follows.  For $0 \leq i, j \leq d$, their $(i,j)$-entries are 
$$
P_{ij} = v_j(\theta_{t+i}), \qquad \qquad P^{*}_{ij}  =  v^{*}_j(\theta^{*}_{r+i}),
$$
where $v_j= v_j^{W}, \ v^*_j= v_j^{*W}$ are from Definition \ref{defdefvi}.
\end{definition}

\begin{theorem}
With reference to Definition \ref{pmatrix}, $P^{*}P = \nu(W)I$, where $\nu(W)$ is from Definition \ref{defnu}.
\end{theorem}

\proof We compute the $(i,j)$-entry of $P^{*}P$ using Definition \ref{pmatrix} and (\ref{AWDvi}). We find that this is equal to $(k_i(W))^{-1}$ times the left hand side of (\ref{ortho1}).  Using (\ref{ortho1}), we obtain $P^{*}P = \nu(W)I$. \hfill $\Box$

\begin{theorem}
Let $\flat$ and $\sharp$ be the maps in Definition \ref{maps}.  With reference to Definition \ref{pmatrix}, 
$Y^{\sharp}P = PY^{\flat}$ for $Y \in End(W)$.
\end{theorem}

\proof
By Lemma \ref{transmatrix}, the transition matrix from a standard basis to a dual standard basis for $W$ is a scalar multiple of $P$.  Therefore, $Y^{\sharp}P = PY^{\flat}$. \hfill $\Box$

\section {Two more bases for $W$}
Let $W$ be as in Assumption \ref{W}.  In Sections \ref{sec:2bases} and \ref{seconbici}, we found two bases for $W$ with respect to which $A$ and $A^*$ are represented by tridiagonal and diagonal matrices.  In this section, we will look at two more bases for $W$ with respect to which $A$ and $A^*$ are represented by lower bidiagonal and upper bidiagonal matrices. 

\begin{definition} \rm\label{taueta}
For $0 \leq i \leq d$, define $\tau_i = \tau_i^{W}\!, \ \tau^*_i = \tau_i^{*W}\!,\ \eta_i = \eta_i^{W}\!, \ \eta^*_i = \eta_i^{W}$ in $\mathbb{C}[\lambda]$ as follows: 
\begin{align*}
\tau_{i} &= {\displaystyle \prod_{h = 0}^{i-1}(\lambda - \theta_{t+h})}, & & \tau^{*}_{i} = \prod_{h =0}^{i-1}(\lambda - \theta^{*}_{r+h}),\\
\eta_i^{} & = {\displaystyle \prod_{h = 0}^{i-1}(\lambda - \theta_{t+d-h})}, & & \eta^{*}_i =  \prod_{h =0}^{i-1}(\lambda - \theta^{*}_{r+d-h}). 
\end{align*} 
Observe that each of $\tau_i,\ \tau^*_i,\ \eta_i,\ \eta^*_i$ is monic of degree $i$.
\end{definition}

\begin{lemma}\label{triviatau}
For $0 \leq i,j \leq d$, 
\begin{enumerate}
\item[\rm(i)] each of $\tau_i(\theta_{t+j}), \tau^*_i(\theta^*_{r+j})$ is $0$ if $j <i$ and nonzero if $j=i$;
\item[\rm(ii)] each of $\eta_i(\theta_{t+j}), \eta^*_i(\theta^*_{r+j})$ is $0$ if $j > d-i$ and nonzero if $j = d-i$.
\end{enumerate}
\end{lemma}

\proof
Immediate from Definition \ref{taueta}. \hfill $\Box$

\begin{lemma}\label {tauibasis}
	Let $v$ be a nonzero vector in $E^*_rW$.  Then $\{\tau_i(A)v \}_{i=0}^d$ is a basis for $W$. 
\end{lemma}

\proof
By Theorem \ref{basis} and Lemma \ref{ppush}, $\{p_i(A)v\}_{i=0}^d$ is a basis for $W$. For $0 \leq i \leq d$, each of $\tau_i$ and $p_i$ is a polynomial of degree $i$.  The result follows. \hfill $\Box$

\begin{definition}\label{Uidef} \rm
For $0 \leq i \leq d$, define
$$U_i = \tau_i(A)E^*_{r}W.$$
For notational convenience, define $U_{-1} = 0$ and $U_{d+1} = 0$.
\end{definition}

\begin{lemma}\label{Uidirect}
With reference to Definition \ref{Uidef}, $U_i$ has dimension $1$ for $0 \leq i \leq d$.  Moreover,
\begin{equation}
W = {\displaystyle \sum_{i=0}^d U_i} \qquad \mbox{\rm(direct sum)}. \label{splitsum}
\end{equation}
\end{lemma}

\proof
Immediate from Lemma \ref{tauibasis} and Definition \ref{Uidef}. \hfill $\Box$

\begin{lemma}\label{basistau}
For $0 \leq i \leq d$,
\begin{enumerate}
\item[\rm(i)]	$\sum_{h=0}^i U_h  =  \sum_{h=0}^i E^{*}_{r+h}W\!$,	
\item[\rm(ii)]	$\sum_{h=i}^d U_h  = \sum_{h=i}^d E_{t+h}W.$	
\end{enumerate}
\end{lemma}

\proof
Let $v$ be a nonzero vector in $E^*_rW$. \\ 
(i)  By Lemma \ref{lemma:ds1}(i),	$\tau_j(A)v$ is contained in $\sum_{h=0}^i E^{*}_{r+h}W$ for $0 \leq j \leq i.$	Hence, $\sum_{h=0}^i U_h \subseteq \sum_{h=0}^i E^{*}_{r+h}W$.  In this inclusion, equality holds since each side has dimension $i+1$.   \\
(ii) For $i \leq j \leq d$, 
\begin{align*}
	\tau_j(A)v \ &= \ \sum_{l=0}^D E_{l}\tau_j(A)v\\
	           \ &= \ \sum_{h=0}^dE_{t+h}\tau_j(A)v\\
	           \ &= \ \sum_{h=0}^d\tau_j(\theta_{t+h})E_{t+h}v\\
  	         \ &= \ \sum_{h=j}^d\tau_j(\theta_{t+h})E_{t+h}v \qquad \qquad \mbox{by Lemma \ref{triviatau}}.
\end{align*}
Hence $\tau_j(A)v \in \sum_{h=i}^d E_{t+h}W$ for $i \leq j \leq d$.  Thus, $\sum_{h=i}^d U_h \subseteq \sum_{h=i}^d E_{t+h}W$.  In this inclusion, equality holds since each side has dimension $i+1$.  \hfill $\Box$

\begin{lemma}\label{Uiother}
For $0 \leq i \leq d$,
$$
U_i = \left(\sum_{h=0}^i E^*_{r+h}W\right) \cap \left(\sum_{h=i}^d E_{t+h}W\right).  
$$
\end{lemma}

\proof By Lemma \ref{Uidirect}, $U_i = (U_0 + U_1 + \cdots + U_i) \cap (U_i + U_{i+1} + \cdots + U_d )$.  Combining this with Lemma \ref{basistau}, we obtain the desired result. \hfill $\Box$

\begin{lemma} \label{raiselower}
For $0 \leq i \leq d$,
\begin{enumerate}
\item[\rm(i)] $(A-\theta_{t+i}I)U_i = U_{i+1}$,
\item[\rm(ii)] $(A^* - \theta^*_{r+i}I)U_i = U_{i-1}$.
\end{enumerate}
\end{lemma}

\proof
(i) Immediate from Definition \ref{Uidef}.\\
(ii) Assume $1 \leq i \leq d$, otherwise, we are done since $U_0 = E^*_rW$.  Let $v$ be a nonzero vector in $E^*_rW$. Since $A^*E^*_{r+i} = \theta^*_{r+i}E^*_{r+i}$, we have $(A^*-\theta^*_{r+i}I)(\sum_{h=0}^iE^*_{r+h}W) \subseteq \sum_{h=0}^{i-1} E^{*}_{r+h}W\!.$  By Lemma \ref{lemma:ds2}(i), we have $(A^*-\theta^*_{r+i}I)(\sum_{h=i}^d E_{t+h}) \subseteq \sum_{h=i-1}^d E_{t+h}W$.  Combining these comments with Lemma \ref{Uiother}, we find that $(A^*-\theta^*_{r+i}I)U_i \subseteq U_{i-1}$. We now show equality holds. Suppose that $(A^*-\theta^*_{r+i}I)U_i \subsetneq U_{i-1}$.  Then  $(A^*-\theta^*_{r+i}I)U_i = 0$ since $\dim U_{i-1} =1$.   Let $W' = U_i + U_{i+1} + \cdots + U_d$.  Observe that $W'$ is nonzero.  By (i), $AW' \subseteq W'$.  Since $(A^{*} - \theta^{*}_{r+i}I)U_i = 0$ and $(A^*-\theta^*_{r+j}I)U_j \subseteq U_{j-1}$ for $i+1 \leq j \leq d$, we find that $A^*W' \subseteq W'.$  Hence $W'$ is a nonzero $T$-submodule of $W$. Since the $T$-module $W$ is irreducible,  $W' = W$\!.  This contradicts (\ref{splitsum}) since $i > 0$. Therefore, $(A^*- \theta^*_{r+i}I)U_i = U_{i-1}$.\hfill $\Box$\\

\indent
By Lemma \ref{raiselower}, for $1 \leq i \leq d$, $U_i$ is invariant under $(A- \theta_{t+i-1}I)(A^{*} - \theta^{*}_{r+i}I)$ and the corresponding eigenvalue is nonzero.  	
	
\begin{definition} \rm\label{def:phii}
For $1 \leq i \leq d$, let $\varphi_i= \varphi_i(W)$ be the eigenvalue of $(A- \theta_{t+i-1}I)(A^{*} - \theta^{*}_{r+i}I)$ corresponding to $U_i$. Observe that $\varphi_i \neq 0$.  We refer to the sequence $\{\varphi_i \}_{i=1}^d$ as the \textit{first split sequence of} $W$\!.  For notational convenience, define $\varphi_0 = 0$.   
\end{definition}

\begin{theorem}\label{wrtsplit1}
With respect to the basis for $W$ in Lemma \ref{tauibasis}, the matrices representing $A$, $A^*$ are
$$\left(
    \begin{array}{cccccc}
    \theta_{t} &              &              &        &           &  {\bf 0} \\
    1          & \theta_{t+1} &              &        &           &          \\
               &  1           & \theta_{t+2} &        &           &          \\
               &              &      \cdot   &  \cdot &           &          \\  
               &              &              &     \cdot   & \theta_{t+d-1}&  \\
    {\bf 0}    &              &              &        & 1         & \theta_{t+d}  \\
    \end{array}
\right), \qquad  \left(
    \begin{array}{cccccc}
    \theta^{*}_{r} &   \varphi_1         &                  &           &                   &  {\bf 0} \\
                   &   \theta^{*}_{r+1}  & \varphi_2        &           &                   &          \\
                   &                     & \theta^{*}_{r+2} & \cdot    &                   &          \\
                   &                     &                  &  \cdot   & \cdot                  &          \\
                   &                     &                  &           & \theta^{*}_{r+d-1}& \varphi_d \\
    {\bf 0}        &                     &                  &           &                   & \theta^{*}_{r+d}  \\
    \end{array}
\right).
$$
\end{theorem}

\proof
Immediate from Definitions \ref{Uidef}, \ref{def:phii} and Lemma \ref{raiselower}. \hfill $\Box$\\

In Lemmas \ref{tauibasis}--\ref{raiselower} and Theorem \ref{wrtsplit1}, we replace $E_{t+i}$ with $E_{t+d-i}$ for $0 \leq i \leq d$ and we routinely obtain the following results.

\begin{lemma}\label {etaibasis}
	Let $v$ be a nonzero vector in $E^*_rW$.  Then $\{\eta_i(A)v \}_{i=0}^d$ is a basis for $W$. 
\end{lemma}

\begin{definition}\label{Uidowndef} \rm
For $0 \leq i \leq d$, define
$$U_i^{\Downarrow} = \eta_i(A)E^*_{r}W.$$
For notational convenience, define $U^{\Downarrow}_{-1} = 0$ and $U^{\Downarrow}_{d+1} = 0$.
\end{definition}	

\begin{lemma}
With reference to Definition \ref{Uidowndef},
\begin{equation}
W = {\displaystyle \sum_{i=0}^d U_i^{\Downarrow}} \qquad \mbox{\rm(direct sum)}. \label{splitsumdown}
\end{equation}
\end{lemma}

\begin{lemma}\label{basiseta}
For $0 \leq i \leq d$,
\begin{enumerate}
\item[\rm(i)]	$\sum_{h= 0}^i U^{\Downarrow}_h =  \sum_{h= 0}^i E^{*}_{r+h}W\!$,  
\item[\rm(ii)]	$\sum_{h= i}^dU^{\Downarrow}_h =  \sum_{h= 0}^{d-i} E_{t+h}W.$
\end{enumerate}
\end{lemma}

\begin{lemma} \label{raiselower2}
For $0 \leq i \leq d$,
\begin{enumerate}
\item[\rm(i)] $(A-\theta_{t+d-i}I)U^{\Downarrow}_i = U^{\Downarrow}_{i+1}$,
\item[\rm(ii)] $(A^* - \theta^*_{r+i}I)U^{\Downarrow}_i = U^{\Downarrow}_{i-1}$.
\end{enumerate}
\end{lemma}

By Lemma \ref{raiselower2}, for $1 \leq i \leq d$, $U^{\Downarrow}_i$ is invariant under $(A- \theta_{t+d-i+1}I)(A^{*} - \theta^{*}_{r+i}I)$ and the corresponding eigenvalue is nonzero. 

\begin{definition} \rm\label{def:phi}
For $1 \leq i \leq d$, let $\phi_i= \phi_i(W)$ be the eigenvalue of $(A- \theta_{t+d-i+1}I)(A^{*} - \theta^{*}_{r+i}I)$ corresponding to $U^{\Downarrow}_i$.  Observe that $\phi_i \neq 0$. We refer to the sequence $\{\phi_i \}_{i=1}^d$ as the \textit{second split sequence of} $W$\!.
\end{definition}

\begin{theorem}\label{wrtsplit2}
With respect to the basis for $W$ in Lemma \ref{etaibasis}, the matrices representing $A,A^*$ are
$$\left(
    \begin{array}{cccccc}
    \theta_{t+d} &              &              &        &           &  {\bf 0} \\
    1          & \theta_{t+d-1} &              &        &           &          \\
               &  1           & \theta_{t+d-2}        &        &           &          \\
               &              &      \cdot   &  \cdot &           &          \\  
               &              &              &     \cdot   & \theta_{t+1}&  \\
    {\bf 0}    &              &              &        & 1         & \theta_{t}  \\
    \end{array}
\right), \qquad \left(
    \begin{array}{cccccc}
    \theta^{*}_{r} &   \phi_1         &                  &           &                   &  {\bf 0} \\
                   &   \theta^{*}_{r+1}  & \phi_2        &           &                   &          \\
                   &                     & \theta^*_{r+2}       & \cdot    &                   &          \\
                   &                     &                  &  \cdot   & \cdot                  &          \\
                   &                     &                  &           & \theta^{*}_{r+d-1}& \phi_d \\
    {\bf 0}        &                     &                  &           &                   & \theta^{*}_{r+d}  \\
    \end{array}
\right).
$$
\end{theorem}
\indent
In \cite[Lemma 12.7]{split}, it was shown that $\{\varphi_i\}_{i=1}^d$ and $\{\phi_i\}_{i=1}^d$ are related by the following:
\begin{eqnarray}
\varphi_i & = & \phi_1 {\displaystyle \sum_{h=0}^{i-1} \frac{\theta_{t+h} - \theta_{t+d-h}}{\theta_{t} - \theta_{t+d}} } + (\theta^{*}_{r+i} - \theta^{*}_r)(\theta_{t+i-1} - \theta_{t+d}) \quad \quad (1 \leq i \leq d), \label{varphi}\\
\phi_i & = & \varphi_1 {\displaystyle \sum_{h=0}^{i-1} \frac{\theta_{t+h} - \theta_{t+d-h}}{\theta_{t} - \theta_{t+d}} } + (\theta^{*}_{r+i} - \theta^{*}_r)(\theta_{t+d-i+1} - \theta_{t}) \quad \quad (1 \leq i \leq d). \label{phi}
\end{eqnarray}
\begin{definition} \rm
By the \textit{parameter array of} $W$, we mean the sequence of scalars
$$(\{\theta_{t+i}\}_{i=0}^d, \{\theta^{*}_{r+i}\}_{i=0}^d, \{\varphi_i\}_{i=1}^d, \{\phi_i\}_{i=1}^d ),$$ where $r, t, d$ are from Assumption \ref{W}, and the $\varphi_i,\ \phi_i$ are from Definitions \ref{def:phii}, \ref{def:phi}.  
\end{definition}

\section {Describing $W$ in terms of its parameter array} \label{everything}
Let $W$ be as in Assumption \ref{W}.  Up until now, we have associated with $W$ a number of polynomials and parameters.  In this section, we will express all these polynomials and parameters in terms of the parameter array $(\{\theta_{t+i}\}_{i=0}^d, \{\theta^{*}_{r+i}\}_{i=0}^d, \{\varphi_i\}_{i=1}^d, \{\phi_i\}_{i=1}^d )$ of $W$. 
Recall the polynomials $\tau_i, \ \tau^*_i, \ \eta_i, \ \eta^*_i$ from Definition \ref{taueta}.  

\begin{theorem}\label{origin}
For $0 \leq i \leq d$,  
\begin{align}
u_i \ &= \ {\displaystyle \sum_{h= 0}^i \frac{\tau^{*}_h(\theta^{*}_{r+i})}{\varphi_1\varphi_2\cdots \varphi_h}\tau_h}, 
\label{uitaui}\\
u^{*}_i \ &= \ {\displaystyle \sum_{h= 0}^i \frac{\tau_h(\theta_{t+i})}{\varphi_1\varphi_2\cdots \varphi_h}\tau^{*}_h},  \label{uitauidual}
\end{align}
where $u_i = u_i^W\!,$ $u^*_i = u_i^{*W}$ are from Definition \ref{defdefui}.
\end{theorem}

\proof
We first verify (\ref{uitaui}).
Since $u_i$ has degree $i$, there exist complex scalars $\{ \alpha_h\}_{h=0}^i$  such that $u_i =  \sum_{h=0}^i \alpha_h \tau_h$.  By Lemma \ref{triviatau}(i),  $\tau_0(\theta_t) = 1$ and $\tau_i(\theta_t) = 0$ for $1 \leq i \leq d$.  From these comments and since $u_i(\theta_t) = 1$, we have $\alpha_0 = 1$.  Now assume $i \geq 1$, otherwise we are done.  Let $v$ be a nonzero vector in $E^{*}_rW$\!.  By Theorem \ref{vpush} and Lemma \ref{viuiki}, $u_i(A)v \in E^{*}_{r+i}W$\!. Thus,  
\begin{align*}
0\ &= \ (A^{*}-\theta^{*}_{r+i}I)u_i(A)v\\
\ &= \ {\displaystyle \sum_{h=0}^i \alpha_h A^{*}\tau_h(A)v - \theta^{*}_{r+i} \sum_{h=0}^i \alpha_h \tau_h(A)v} \\
 \ &= \ {\displaystyle \sum_{h=0}^i \alpha_h(\theta^*_{r+h} \tau_{h}(A)v + \varphi_{h} \tau_{h-1}(A)v) - \theta^{*}_{r+i}\sum_{h=0}^i \alpha_{h} \tau_{h}(A)v} \qquad \mbox{by Theorem \ref{wrtsplit1}}\\
\ &= \  {\displaystyle \sum_{h=0}^{i-1}(\varphi_{h+1}\alpha_{h+1} + \alpha_h\theta^*_{r+h} - \theta^{*}_{r+i}\alpha_h)\tau_h(A)v}.
\end{align*}
By Lemma \ref{tauibasis}, $\{\tau_h(A)v\}_{h=0}^{i-1}$ are linearly independent. Thus, $\varphi_{h+1}\alpha_{h+1} + \alpha_h\theta^*_{r+h} - \theta^{*}_{r+i}\alpha_h = 0$ for $0 \leq h < i$.  From this recursive equation and the fact that $\alpha_0 = 1$, we find that $\alpha_h = \tau^*_h(\theta^{*}_{r+i})/(\varphi_1\varphi_2 \cdots \varphi_h)$ for $0 \leq h \leq i$.  Therefore, (\ref{uitaui}) holds.  We now prove (\ref{uitauidual}).  
Let $f_i$ be the polynomial on the right in (\ref{uitauidual}).  Using (\ref{uitaui}), we find that $f_i(\theta^*_{r+j}) = u_j(\theta_{t+i})$ for $0 \leq j \leq i$.  
By Theorem \ref{askey}, $u^*_i(\theta^*_{r+j}) = u_j(\theta_{t+i})$.  Therefore, $f_i(\theta^*_{r+j}) = u^*_i(\theta^*_{r+j})$ for $0 \leq j \leq i$. 
By this and since $u^*_i, f_i$ have degree $i$, we find that $u^*_i = f_i$.    \hfill $\Box$

\begin{lemma} \label{sarisari}
For $0 \leq i \leq d$,
\begin{align}
p_i(\theta_t)& = {\displaystyle \frac{\varphi_1\varphi_2\cdots\varphi_i}{\tau^{*}_i(\theta^{*}_{r+i})}}, \qquad \qquad p^{*}_i(\theta^{*}_r) = {\displaystyle \frac{\varphi_1\varphi_2\cdots\varphi_i}{\tau_i(\theta_{t+i})}},  \label{piulit}
\end{align}
where $p_i = p_i^W, \ p^*_i = p_i^{*W}$ are from Definitions \ref{polypi}, \ref{polypi*}.
\end{lemma}

\proof
We first prove the equation on the left in (\ref{piulit}).  We compute the coefficient of $\lambda^i$ in $u_i$ in two ways: one way using (\ref{uitaui}) and another way using Definition \ref{defdefui}.  Comparing the results, we obtain the equation on the left in (\ref{piulit}). Argue similarly to obtain the equation on the right in (\ref{piulit}). \hfill $\Box$

\begin{theorem}\label{thm:biprod}
For $0 \leq i  \leq d-1$,
\begin{equation}
b_i(W) = \varphi_{i+1}{\displaystyle \frac{\tau^{*}_{i}(\theta^{*}_{r+i})}{\tau^{*}_{i+1}(\theta^{*}_{r+i+1})}},  \qquad \qquad b^{*}_i(W) = \varphi_{i+1}{\displaystyle \frac{\tau_{i}(\theta_{t+i})}{\tau_{i+1}(\theta_{t+i+1})}}. \label{biprod}
\end{equation} 
The $b_i(W), \ b^*_i(W)$ are from Definition \ref{def:bi&ci}.  
\end{theorem}

\proof
Immediate from Theorem \ref{formbici}(i), (iii) and Lemma \ref{sarisari}. \hfill $\Box$

\begin{theorem}\label{thm:ai}
With reference to Definition \ref{defxiai},
\begin{align}
a_0(W) \ &=\ {\displaystyle\theta_{t} + \frac{\varphi_{1}}{\theta^{*}_{r} - \theta^{*}_{r+1}}},
& a_d(W) \ &= \ {\displaystyle\theta_{t+d} + \frac{\varphi_d}{\theta^{*}_{r+d} - \theta^{*}_{r+d-1}}}, \label{a01}\\
a_0(W) \ &= \  {\displaystyle\theta_{t+d} +  \frac{\phi_{1}}{\theta^{*}_{r} - \theta^{*}_{r+1}}},
& a_d(W) \ &= \ {\displaystyle\theta_{t} + \frac{\phi_d}{\theta^{*}_{r+d} - \theta^{*}_{r+d-1}}}.\label{a02} 
\end{align}
For $1 \leq i \leq d-1$,
\begin{align}
a_i(W)\ &=\ {\displaystyle\theta_{t+i} + \frac{\varphi_i}{\theta^{*}_{r+i} - \theta^{*}_{r+i-1}} + \frac{\varphi_{i+1}}{\theta^{*}_{r+i} - \theta^{*}_{r+i+1}}} \label{ai1} \\
\ &=\ {\displaystyle\theta_{t+d-i} + \frac{\phi_i}{\theta^{*}_{r+i} - \theta^{*}_{r+i-1}} + \frac{\phi_{i+1}}{\theta^{*}_{r+i} - \theta^{*}_{r+i+1}}}.\label{ai2}
\end{align}
\end{theorem}

%

\proof
To obtain (\ref{ai1}), we compute the coefficient of $\lambda^i$ in $u_{i+1}$ in two ways.  One way is using Lemma \ref{bipi} and Lemma \ref{recurrenceui}.  Using this approach, we find that the coefficient is equal to 
\begin{equation}
-\sum_{l=0}^i\frac{ a_l(W)}{p_{i+1}(\theta_t)} \label{firstcoeff}.
\end{equation}
Another way is using  (\ref{uitaui}).  Using this approach, the coefficient is equal to 
\begin{equation}
{\displaystyle\frac{\tau^*_{i}(\theta^*_{r+i+1})}{\varphi_1\varphi_2 \cdots \varphi_i} 
-  \sum_{l=0}^i \theta_{t+l}\frac{\tau^*_{i+1}(\theta^*_{r+i+1})}{\varphi_1\varphi_2 \cdots \varphi_{i+1}}}.
\label{2ndcoeff}
\end{equation}
Evaluating (\ref{firstcoeff}) using (\ref{piulit}) and comparing the result with (\ref{2ndcoeff}), we obtain (\ref{ai1}).    Similarly, we obtain the two equations in (\ref{a01}).  We now prove (\ref{ai2}).  Observe that by Definitions \ref{defxiai} and \ref{def:phi}, replacing $E_{t+i}$ with $E_{t+d-i}$ for $0 \leq i \leq d$ has the effect of switching $(a_i(W), \theta_{t+i}, \varphi_i)$ to $(a_i(W), \theta_{t+d-i}, \phi_i)$. Applying this switching to (\ref{ai1}), we obtain (\ref{ai2}).  Similarly, we obtain the two equations in (\ref{a02}). \hfill $\Box$
\begin{theorem}
With reference to Definition \ref{defxiai},
\begin{align}
a^{*}_0(W) \ &= \ {\displaystyle \theta^{*}_{r} +  \frac{\varphi_{1}}{\theta_{t} - \theta_{t+1}}},
& a^{*}_d(W) \ &= \ {\displaystyle \theta^{*}_{r+d} + \frac{\varphi_d}{\theta_{t+d} - \theta_{t+d-1}}} \label{a01*}, \\
a^*_0(W)\ & =\ {\displaystyle \theta^{*}_{r+d} +  \frac{\phi_{d}}{\theta_{t} - \theta_{t+1}}}, 
& a^{*}_d(W) \ & =\ {\displaystyle \theta^{*}_{r} + \frac{\phi_1}{\theta_{t+d} - \theta_{t+d-1}}}.
\label{a02*}
\end{align}
For $1 \leq i \leq d-1$,
\begin{align}
a^{*}_i(W)\ &= \ {\displaystyle \theta^{*}_{r+i} + \frac{\varphi_i}{\theta_{t+i} - \theta_{t+i-1}} + \frac{\varphi_{i+1}}{\theta_{t+i} - \theta_{t+i+1}}} \label{ai1*}\\
\ & =\ {\displaystyle \theta^{*}_{r+d-i} + \frac{\phi_{d-i+1}}{\theta_{t+i} - \theta_{t+i-1}} + \frac{\phi_{d-i}}{\theta_{t+i} - \theta_{t+i+1}}}. \label{ai2*}
\end{align}
\end{theorem}

\proof
To obtain (\ref{a01*}) and (\ref{ai1*}) argue similarly as in the proof of (\ref{ai1}).  We now prove (\ref{ai2*}).  By Definitions \ref{defxiai} and \ref{def:phi}, replacing $E_{t+i}$ with $E_{t+d-i}$ for $0 \leq i \leq d$ has the effect of switching $(a^*_{i}(W), \theta_{t+i}, \varphi_i)$ to $(a^*_{d-i}(W), \theta_{t+d-i}, \phi_i)$.  Applying this switching to (\ref{ai1*}), we obtain  
\begin{equation}
a^{*}_{d-i}(W) = {\displaystyle \theta^{*}_{r+i} + \frac{\phi_i}{\theta_{t+d-i} - \theta_{t+d-i+1}} + \frac{\phi_{i+1}}{\theta_{t+d-i} - \theta_{t+d-i-1}}}. \label {ai3*}
\end{equation}
Changing $i$ to $d-i$ in (\ref{ai3*}), we obtain (\ref{ai2*}). \hfill $\Box$

\begin{theorem} \label{formofvarphi}
For $1 \leq i \leq d$, $\varphi_i$ is equal to each of the following:
\begin{align}
& (\theta^{*}_{r+i} - \theta^{*}_{r+i-1}) {\displaystyle\sum_{j=0}^{i-1} (\theta_{t+j} -a_j(W))}, &  & (\theta^{*}_{r+i-1} - \theta^*_{r+i}){\displaystyle \sum_{j=i}^d (\theta_{t+j} -a_j(W))}, \label{varphi1}\\
&  (\theta_{t+i} - \theta_{t+i-1}) {\displaystyle \sum_{j=0}^{i-1} (\theta^{*}_{r+j} -a^{*}_j(W))},   &  & (\theta_{t+i-1} - \theta_{t+i}){\displaystyle \sum_{j=i}^d (\theta^{*}_{r+j} -a^{*}_j(W))}.  \label{varphi2}
\end{align} 
The $a_h(W), \ a^*_h(W)$ are from Definition \ref{defxiai}.
\end{theorem}

\proof
To obtain the expression on the left in (\ref{varphi1}), solve for $\varphi_{i}$ recursively using (\ref{ai1}).  From this and Lemma \ref{sumofai}(i), we obtain the expression on the right in (\ref{varphi1}).  The remaining assertions can be similarly shown.  \hfill $\Box$

\begin{theorem}\label{formofphi}
For $1 \leq i \leq d$, $\phi_i$ is equal to each of the following:
\begin{align}
& (\theta^{*}_{r+i} - \theta^{*}_{r+i-1}) {\displaystyle\sum_{j=0}^{i-1} (\theta_{t+d-j} -a_j(W))}, &  & (\theta^{*}_{r+i-1} - \theta^*_{r+i}){\displaystyle \sum_{j=i}^d (\theta_{t+d-j} -a_j(W))}, \label{phi1} \\
&  (\theta_{t+d-i} - \theta_{t+d-i+1}) {\displaystyle \sum_{j=0}^{i-1} (\theta^{*}_{r+j} -a^{*}_{d-j}(W))},   &  & (\theta_{t+d-i+1} - \theta_{t+d-i}){\displaystyle \sum_{j=i}^d (\theta^{*}_{r+j} -a^{*}_{d-j}(W))}. \label{phi2}
\end{align} 
The $a_h(W), \ a^*_h(W)$ are from Definition \ref{defxiai}.
\end{theorem}

\proof
Similar to the proof of Theorem \ref{formofvarphi}. \hfill $\Box$ 

\begin{theorem}

For $0 \leq i \leq d$, the polynomial $p_i = p_i^W$ from Definition \ref{polypi} is equal to both
\begin{equation}
{\displaystyle \sum_{h=0}^i \frac{\varphi_1\varphi_2\cdots \varphi_i\tau^{*}_h(\theta^{*}_{r+i})}{\varphi_1\varphi_2 \cdots\varphi_h\tau^{*}_i(\theta^{*}_{r+i})}\tau_h}, \qquad \qquad {\displaystyle \sum_{h=0}^i \frac{\phi_1\phi_2\cdots \phi_i\tau^{*}_h(\theta^{*}_{r+i})}{\phi_1\phi_2 \cdots\phi_h\tau^{*}_i(\theta^{*}_{r+i})}\eta_h }. \label{pitaui}
\end{equation}
\end{theorem}

\proof
The expression on the left in (\ref{pitaui}) is equal to $p_i$ by Definition \ref{defdefui}, (\ref{uitaui}), and the equation on the left in (\ref{piulit}). To show that $p_i$ is equal to the expression on the right in (\ref{pitaui}), write $u_i$ as a linear combination of $\{\eta_h \}_{h=0}^i.$  Arguing as in the proof of (\ref{uitaui}), we find that
\begin{equation}
u_i = u_i(\theta_{t+d}){\displaystyle \sum_{h= 0}^i \frac{\tau^{*}_h(\theta^{*}_{r+i})}{\phi_1\phi_2\cdots \phi_h}\eta_h}. \label{anotherforui}
\end{equation}
To find $u_i(\theta_{t+d})$, we compute the coefficient of $\lambda^i$ in $u_i$ in two ways: one way is using (\ref{uitaui}) and another way is using (\ref{anotherforui}).  Comparing these results we obtain
\begin{equation}
u_i(\theta_{t+d}) = \frac{\phi_1\phi_2\cdots \phi_i}{\varphi_1\varphi_2\cdots \varphi_i}. \label{uithetat+d}
\end{equation}
Evaluating $p_i$ using Definition \ref{defdefui}, ({\ref{anotherforui}),  (\ref{uithetat+d}) and the equation on the left in (\ref{piulit}), we find that $p_i$ is equal to the expression on the right in (\ref{pitaui}). \hfill $\Box$  

\begin{theorem} \label{thm:pi*}
For $0 \leq i \leq d$, the polynomial $p^*_i = p_i^{*W}$ from Definition \ref{polypi*} is equal to both
\begin{equation}
{\displaystyle \sum_{h=0}^i \frac{\varphi_1\varphi_2\cdots \varphi_i\tau_h(\theta_{t+i})}{\varphi_1\varphi_2 \cdots\varphi_h\tau_i(\theta_{t+i})}\tau^{*}_h}, \qquad \qquad \label{pitauidual}
{\displaystyle \sum_{h=0}^i \frac{\phi_{d}\phi_{d-1}\cdots \phi_{d-i+1}\tau_h(\theta_{t+i})}{\phi_d\phi_{d-1} \cdots\phi_{d-h+1}\tau_i(\theta_{t+i})}\eta^{*}_h }. 
\end{equation}
\end{theorem}

\proof
The expression on the left in (\ref{pitauidual}) is equal to $p^*_i$ by Definition \ref{defdefui}, (\ref{uitauidual}), and the equation on the right in (\ref{piulit}). We now prove that $p^*_i$ is equal to the expression on the right in (\ref{pitauidual}).  Comparing the equation on the left in (\ref{phi1}) and the equation on the right in (\ref{phi2}), we find that interchanging $A$ and $A^*$ has the effect of switching $\phi_i$ to $\phi_{d-i+1}$ for $1 \leq i \leq d$.  Applying this switching to the sum on the right in (\ref{pitaui}), we obtain the sum on the right in (\ref{pitauidual}).   \hfill $\Box$

\begin{lemma}
For $0 \leq i \leq d$, 
$$
p_i(\theta_{t+d}) = {\displaystyle\frac{\phi_1\phi_2\cdots \phi_i}{\tau^{*}_i(\theta^{*}_{r+i})}}, \qquad \qquad
p^{*}_i(\theta^{*}_{r+d}) = {\displaystyle \frac{\phi_{d}\phi_{d-1}\cdots \phi_{d-i+1}}{\tau_i(\theta_{t+i})}},
$$
where $p_i = p_i^W, \ p^*_i = p_i^{*W}$ are from Definitions \ref{polypi}, \ref{polypi*}.
\end{lemma}

\proof
Immediate from the right side of lines (\ref{pitaui}) and (\ref{pitauidual}). \hfill $\Box$



\begin{theorem} \label{thm:ciprod}
For $1 \leq i \leq d$,
\begin{equation}
c_i(W) =  \phi_{i}{\displaystyle \frac{\eta^{*}_{d-i}(\theta^{*}_{r+i})}{\eta^{*}_{d-i+1}(\theta^{*}_{r+i-1})}}, \qquad \qquad
c^{*}_i(W) = \phi_{d-i+1}{\displaystyle \frac{\eta_{d-i}(\theta_{t+i})}{\eta_{d-i+1}(\theta_{t+i-1})}}. \label{ciprod} 
\end{equation}
The $c_i(W), \ c^*_i(W)$ are from Definition \ref{def:bi&ci}.
\end{theorem}

\proof
We first verify the equation on the right in (\ref{ciprod}).  By (\ref{matrixrepstandard*}), replacing $E_{t+i}$ with $E_{t+d-i}$ for $0 \leq i \leq d$ switches $b^*_i(W)$ and $c^*_{d-i}(W)$.  Applying this switching to the equation on the right in (\ref{biprod}), we find that for $0 \leq i \leq d-1$, 
\begin{equation}
c^*_{d-i}(W) = {\displaystyle \phi_{i+1}\frac{\eta_{i}(\theta_{t+d-i})}{\eta_{i+1}(\theta_{t+d-i-1})}}. \label{eq:1} 
\end{equation}
Changing $i$ to $d-i$ in (\ref{eq:1}), we obtain the equation on the right in (\ref{ciprod}).  We now verify the equation on the left in (\ref{ciprod}).  Recall from the proof of Theorem \ref{thm:pi*} that interchanging $A$ and $A^*$ switches $\phi_i$ and $\phi_{d-i+1}$.  Applying this switching to the equation on the right in (\ref{ciprod}), we obtain the equation on the left in (\ref{ciprod}). \hfill $\Box$

\begin{theorem}
With reference to Definition \ref{defnu},
\begin{equation}
\nu(W)  = \frac{\eta_d(\theta_t)\eta^*_{d}(\theta^*_r)}{\phi_1\phi_2 \cdots \phi_d}. \label{formofnu}
\end{equation}
\end{theorem}

\proof
Let $0 \neq v  \in E^*_rW\!.$  By Theorem \ref{wrtsplit2}, $(A^*- \theta^*_{r+i}I)\eta_i(A)v = \phi_i \eta_{i-1}(A)v$ for $1 \leq i \leq d$.  Hence, $\eta^*_d(A^*)\eta_d(A)v = \phi_1 \phi_2 \cdots \phi_d v$.  By (\ref{eq:explicitEi}) and (\ref{eq:formpidem}), on $W$ we have $\eta_d(A)=\eta_{d}(\theta_t)E_t $ and $\eta^*_d(A^*)=\eta^*_{d}(\theta^*_r)E^*_r$.  Thus, $\eta^*_d(A^*)\eta_d(A)v = \eta^*_d(\theta^*_r)\eta_d(\theta_{t})E^*_rE_tv.$  From these comments and since $v \in E^*_rW$, we obtain $\phi_1 \phi_2 \cdots \phi_dv = \eta^*_d(\theta^*_r)\eta_d(\theta_{t})E^*_rE_tE^*_rv.$  
Evaluate $E^*_rE_tE^*_r$ using Theorem \ref{thm:nu}(ii). The result follows.  \hfill $\Box$

\begin{theorem}
With reference to Definitions \ref{defxiai} and \ref{ki},
\begin{align} 
k_i(W) \ & = \ \frac{\varphi_1\varphi_2\cdots \varphi_i}{\phi_1\phi_2 \cdots \phi_i} \, \frac{\eta^*_d(\theta^*_r)}{\tau^*_i(\theta^*_{r+i})\eta^*_{d-i}(\theta^*_{r+i})}  \ \qquad  \qquad \quad (0 \leq i \leq d), \label{form:ki}\\
k_i^*(W) \ & = \ \frac{\varphi_1\varphi_2\cdots \varphi_i}{\phi_d\phi_{d-1} \cdots \phi_{d-i+1}}\, \frac{\eta_d(\theta_t)}{\tau_i(\theta_{t+i})\eta_{d-i}(\theta_{t+i})}  \ \ \qquad (0 \leq i \leq d), \label{form:ki^*}\\
x_i(W) \ & = \ \varphi_i \phi_i \frac{\tau^*_{i-1}(\theta^*_{r+i-1})\eta^*_{d-i}(\theta^*_{r+i})}{\tau^*_i(\theta^*_{r+i})\eta^*_{d-i+1}(\theta^*_{r+i-1})} \qquad \qquad \qquad (1 \leq i \leq d), \label{form:xi}\\
x^*_i(W) \ & = \   \varphi_i \phi_{d-i+1} \frac{\tau_{i-1}(\theta_{t+i-1})\eta_{d-i}(\theta_{t+i})}{\tau_i(\theta_{t+i})\eta_{d-i+1}(\theta_{t+i-1})}\qquad \qquad \ (1 \leq i \leq d). \label{form:xi^*}
\end{align}
\end{theorem}

\proof
Evaluate the equations in (\ref{kiform}), (\ref{kidualform}) and in Lemma \ref{propbici}(i), (iii), using Theorems \ref{thm:biprod} and \ref{thm:ciprod}.\hfill $\Box$\\

For the rest of this section, we will find alternative formulae for the intersection and dual intersection numbers of $W$.  The reason in doing this is that the formulae given in Theorems \ref{thm:biprod} and \ref{thm:ciprod} involve huge products which may not be easy to compute.  We will need the following lemma.

\begin{lemma} \label{3terms}
For $0 \leq i \leq d$,
\begin{align}
c_i(W)\tau^{*}_1(\theta^{*}_{r+i-1}) + a_i(W)\tau^{*}_1(\theta^{*}_{r+i}) + b_i(W)\tau^{*}_1(\theta^{*}_{r+i+1}) \ &= \ \varphi_1 + \theta_{t+1}\tau^{*}_1(\theta^{*}_{r+i}), \label{3termtau}\\
c^{*}_i(W)\tau_1(\theta_{t+i-1}) + a^{*}_i(W)\tau_1(\theta_{t+i}) + b^{*}_i(W)\tau_1(\theta_{t+i+1}) \ &= \ \varphi_1 + \theta^{*}_{r+1}\tau_1(\theta_{t+i}). \label{3termtaudual}
\end{align}
The $a_i(W), \ b_i(W), \ c_i(W)$ (resp.\! $a^*_i(W), \ b^*_i(W), \ c^*_i(W)$) are the intersection numbers (resp.\! dual intersection numbers) of $W$.
\end{lemma}

\proof
By (\ref{recforui*}), $u^*_1 = (\lambda - a^*_0(W))/b^*_0(W)$.  Use this to evaluate (\ref{difeqn1}) with $j=1$. Eliminate $a^*_0(W)$ in the resulting equation using the expression on the left of (\ref{a01*}). Simplify using Lemma \ref{propbici}(ii) to obtain (\ref{3termtau}).  The proof of (\ref{3termtaudual}) is similar.  \hfill $\Box$

\begin{theorem} \label{formofintarray}
The intersection numbers of $W$ are as follows:
\begin{align} 
b_0(W) &= \frac{\varphi_1}{\theta^{*}_{r+1} - \theta^{*}_r}, \label{formofb0}\\
b_i(W) &= \frac{(\theta_t - a_i(W))(\theta^{*}_{r+i} - \theta^{*}_{r+i-1})+(\theta_{t} - \theta_{t+1})(\theta^{*}_{r} - \theta^{*}_{r+i}) + \varphi_1}{\theta^{*}_{r+i+1} - \theta^{*}_{r+i-1}} \qquad (1 \leq i \leq d-1), \label{formofbi}\\
c_i(W) &= \frac{(\theta_t - a_i(W))(\theta^{*}_{r+i} - \theta^{*}_{r+i+1})+(\theta_{t} - \theta_{t+1})(\theta^{*}_{r} - \theta^{*}_{r+i}) + \varphi_1}{\theta^{*}_{r+i-1} - \theta^{*}_{r+i+1}} \qquad (1 \leq i \leq d-1), \label{formofci}\\
c_d(W) & =  \frac{\varphi_1 + (\theta_{t+1} - \theta_t)(\theta^{*}_{r+d}- \theta^{*}_{r})}{\theta^{*}_{r+d-1} - \theta^{*}_{r+d}}. \label{formofcd}
\end{align}
To obtain $b^*_i(W)$ and $c_i^*(W)$, replace $(\theta_{t+j}, \ \theta^*_{r+j}, \  a_j(W))$ with $(\theta^*_{r+j},\  \theta_{t+j}, \  a^*_j(W))$. 
\end{theorem}

\proof
To obtain (\ref{formofb0}), eliminate $a_0(W)$ in the equation on the left of (\ref{a01}) using Lemma \ref{propbici}(ii). 
To obtain (\ref{formofbi}) and (\ref{formofci}), solve the system of equations in Lemmas \ref{propbici}(ii) and (\ref{3termtau}).  To obtain (\ref{formofcd}), set $i=d$ in (\ref{3termtau}) and eliminate $a_d(W)$  using Lemma \ref{propbici}(ii).  The proof of the assertion regarding the dual intersection numbers of $W$ is similar.  \hfill $\Box$\\

\indent
By Theorems \ref{thm:ai}, \ref{formofintarray}, the intersection numbers (resp.\! dual intersection numbers) of $W$ can be expressed in terms of the parameter array of $W$.  By (\ref{varphi}) and (\ref{phi}), the parameter array of $W$ is determined by the eigenvalue sequence of $W$, dual eigenvalue sequence of $W$, and $\varphi_1(W)$.  Hence, we now solve for the intersection numbers (resp.\! dual intersection numbers) of $W$ in terms of these parameters.  But first we need the following lemmas.     

\begin{lemma} \label{prelim}
Assume $d \geq 2$.  Then the scalar $\varphi_2$ is equal to both
\begin{equation}
\varphi_1(1\! +\!  \frac{\theta_{t+1}\!  -\!  \theta_{t+d-1}}{\theta_{t}\!  -\!  \theta_{t+d}})\!  +\!  (\theta^*_{r+1}\!  -\!  \theta^*_{r})(\theta_{t+d}\! +\! \theta_{t+d-1}\! -\! \theta_{t}\! -\! \theta_{t+1})\!  +\!  (\theta^*_{r+2}\!  -\!  \theta^*_{r})(\theta_{t+1} \! - \theta_{t+d}), \label{anotherformofvarphi2}
\end{equation}
\begin{equation}
\varphi_1(1\! +\! \frac{\theta^*_{r+1}\! -\! \theta^*_{r+d-1}}{\theta^*_{r}\! -\! \theta^*_{r+d}})\! +\! (\theta_{t+1}\! -\! \theta_{t})(\theta^*_{r+d}\!+\!\theta^*_{r+d-1}\!-\!\theta^*_{r}\!-\!\theta^*_{r+1})\! +\! (\theta_{t+2}\! -\! \theta_{t})(\theta^*_{r+1}\! -\! \theta^*_{r+d}). \label{formofvarphi2}
\end{equation}
\end{lemma}

\proof
To obtain (\ref{anotherformofvarphi2}), set $i=2$ in (\ref{varphi}) and evaluate $\phi_1$ using (\ref{phi}). Comparing the formula for $\varphi_i$ on the left in lines (\ref{varphi1}) and (\ref{varphi2}), we find that interchanging $A$ and $A^*$ has no effect on $\varphi_i$ for $1 \leq i \leq d$.  Applying this switching to (\ref{anotherformofvarphi2}), we obtain (\ref{formofvarphi2}).   \hfill $\Box$

\begin{lemma}
Assume $d \geq 2$. Then for $0 \leq i \leq d$,
\begin{align}
c_i(W)\tau^{*}_2(\theta^{*}_{r+i-1}) + a_i(W)\tau^{*}_2(\theta^{*}_{r+i}) + b_i(W)\tau^{*}_2(\theta^{*}_{r+i+1}) \ &= \ \varphi_2\tau^{*}_{1}(\theta^{*}_{r+i}) + \theta_{t+2}\tau^{*}_2(\theta^{*}_{r+i}), \label{3termtau2}\\
c^{*}_i(W)\tau_2(\theta_{t+i-1}) + a^{*}_i(W)\tau_2(\theta_{t+i}) + b^{*}_i(W)\tau_2(\theta_{t+i+1}) \ &= \ \varphi_2\tau_{1}(\theta_{t+i}) + \theta^{*}_{r+2}\tau_2(\theta_{t+i}).\label{3termtaudual2}
\end{align}
The $a_i(W), \ b_i(W), \ c_i(W)$ (resp.\! $a^*_i(W), \ b^*_i(W), \ c^*_i(W)$) are the intersection numbers (resp.\! dual intersection numbers) of $W$.
\end{lemma}

\proof
Eliminating $u^*_2$ in (\ref{difeqn1}) with $j=2$ using (\ref{uitauidual}), we obtain
\begin{align}
c_i(W) &+ a_i(W) +b_i(W) + \frac{\tau_1(\theta_{t+2})}{\varphi_1}(c_i(W)\tau^{*}_1(\theta^{*}_{r+i-1}) + a_i(W)\tau^{*}_1(\theta^{*}_{r+i}) + b_i(W)\tau^{*}_1(\theta^{*}_{r+i+1})) \label{above}\\
&+\frac{\tau_2(\theta_{t+2})}{\varphi_1\varphi_2}(c_i(W)\tau^{*}_2(\theta^{*}_{r+i-1}) + a_i(W)\tau^{*}_2(\theta^{*}_{r+i}) + b_i(W)\tau^{*}_2(\theta^{*}_{r+i+1})) \nonumber\\
&= \theta_{t+2}(1+\frac{\tau_1(\theta_{t+2})}{\varphi_1}\tau^*_1(\theta^*_{r+i}) + \frac{\tau_2(\theta_{t+2})}{\varphi_1\varphi_2}\tau^*_2(\theta^*_{r+i}) ).  \nonumber
\end{align}
Simplify the first three terms of (\ref{above}) using Lemma \ref{propbici}(ii).  Evaluating the coefficient of $\tau_1(\theta_{t+2})/\varphi_1$ in (\ref{above}) using (\ref{3termtau}), we routinely obtain (\ref{3termtau2}).  The proof of (\ref{3termtaudual2}) is similar.  \hfill $\Box$

\begin{theorem} \label{formofintarray2}
The intersection numbers of $W$ are as follows:
\begin{align}
b_0(W) \ & = \ \frac{\varphi_1}{\theta^{*}_{r+1} - \theta^{*}_r}, \label{boagain}\\
b_i(W)  \ & = \ {\displaystyle\frac{\varphi_1f^{+}_i +  g^{+}_i}{(\theta^{*}_{r+i+1} - \theta^{*}_{r+i})(\theta^{*}_{r+i+1}-\theta^{*}_{r+i-1} )}}, \quad \quad (1 \leq i \leq d-1), \label{biform}\\ 
c_i(W)  \ & = \ {\displaystyle\frac{\varphi_1f^{-}_i + g^{-}_i}{(\theta^{*}_{r+i-1} - \theta^{*}_{r+i})(\theta^{*}_{r+i-1}-\theta^{*}_{r+i+1} )}} \ \quad \quad (1 \leq i \leq d-1), \label{ciform}\\
c_d(W)  \ & = \ \frac{\varphi_1 + (\theta_{t+1}-\theta_{t})(\theta^*_{r+d} - \theta^*_{r})}{\theta^*_{r+d-1} - \theta^*_{r+d}}, \label{cdagain}
\end{align}
where
\begin{align*}
f^{\pm}_i \ & = \ \theta^{*}_{r+1} - \theta^{*}_{r+i\mp 1} - \frac{(\theta^{*}_{r+i} - \theta^{*}_{r})(\theta^{*}_{r+1}- \theta^{*}_{r+d-1})}{(\theta^{*}_{r+d} - \theta^{*}_{r})}, \\
g^{\pm}_i \ & = \ (\theta^{*}_{r+i} - \theta^{*}_r)((\theta_{t+2} - \theta_{t+1})(\theta^{*}_{r+i}- \theta^{*}_{r+d}) - (\theta_{t+1}- \theta_t)(\theta^{*}_{r+i \mp 1} - \theta^{*}_{r+d-1})),
\end{align*}
provided $d \geq 2$.
To obtain $b^*_i(W)$ and $c^*_i(W)$, replace $(\theta_{t+j}, \theta^*_{r+j})$ with $(\theta^*_{r+j}, \theta_{t+j})$.
\end{theorem}

\proof
 Observe that (\ref{boagain}), (\ref{cdagain}) are  (\ref{formofb0}), (\ref{formofcd}). To obtain (\ref{biform}) and (\ref{ciform}), eliminate $a_i(W)$ in (\ref{3termtau}) and (\ref{3termtau2}) using Lemma \ref{propbici}(ii).  Then for $1 \leq i \leq d-1$, 
\begin{equation}
c_i(W)(\theta^*_{r+i-1} - \theta^*_{r+i}) + b_i(W)(\theta^*_{r+i+1} - \theta^*_{r+i}) = \varphi_1 + (\theta_{t+1} - \theta_{t})(\theta^*_{r+i} - \theta^*_r), \label{sys2}
\end{equation}
\begin{equation}
c_i(W)h_i(\theta^*_{r+i-1}) + b_i(W)h_i(\theta^*_{r+i+1}) =\varphi_2(\theta^*_{r+i} -\theta^*_{r}) +  (\theta_{t+2} - \theta_t)(\theta^*_{r+i}- \theta^*_{r})(\theta^*_{r+i}- \theta^*_{r+1}), \label{eqhi}
\end{equation}
where
$$h_i(\lambda) = (\lambda - \theta^*_{r+i})(\lambda + \theta^*_{r+i} - \theta^*_{r+1} - \theta^*_{r}).$$ 
 Eliminate $\varphi_2$ in (\ref{eqhi}) using (\ref{formofvarphi2}). Solving the system of equations (\ref{sys2}), (\ref{eqhi}), we obtain (\ref{biform}) and (\ref{ciform}).
Argue similarly and evaluate $\varphi_2$ using (\ref{anotherformofvarphi2}) to obtain the formula for the dual intersection numbers of $W$.  \hfill $\Box$

\begin{lemma} \label{indep}
Given vertices $y,z$ in $X$, let $W$ (resp. $W'$) be the trivial $T(y)$-module (resp. $T(z)$-module) of $\Gamma$.  Then $W$ and $W'$ have the same parameter array. \hfill $\Box$
\end{lemma}

\proof
By Lemma \ref{fortrivial}, $W$ and $W'$ have the same intersection numbers and dual intersection numbers. Thus, $W$ and $W'$ both have eigenvalue sequence $\{\theta_i\}_{i=0}^D$ and dual eigenvalue sequence $\{\theta^*_i\}_{i=0}^D$.  By (\ref{boagain}), $\varphi_1(W) = \varphi_1(W')$. Using (\ref{varphi}) and (\ref{phi}), we find that $W$ and $W'$ have the same first split sequence and second split sequence.  
\hfill $\Box$

\begin{definition}\label{PAgamma} \rm
By the \textit{parameter array of} $\Gamma$, we mean the parameter array of the trivial $T(x)$-module. Observe that this parameter array is independent of the choice of $x$ by Lemma \ref{indep}.
\end{definition}

\section {Isomorphism Classes of Thin Irreducible T-modules}
In Corollary \ref{isoclass}, we mentioned some set of scalars needed to determine the isomorphism class of a thin irreducible $T\!$-module.  As we have seen in Theorem \ref{formofintarray2}, there are many relations among these scalars.We now consider a much smaller set of scalars needed to determine the isomorphism class.  Let us first consider some equations from (\ref{varphi}), (\ref{phi}).  

\begin{lemma} \label{mgaphi}
Let $W$ be as in Assumption \ref{W}.  Let $\{\varphi_i\}_{i=1}^d, \{\phi_i\}_{i=1}^d$ denote the first split sequence and second split  sequence of $W$\!, respectively.  Then 
\begin{align*}
\phi_1 \ &= \ \varphi_1 +(\theta^*_{r+1}-\theta^*_r)(\theta_{t+d} - \theta_t),\\
\phi_d \ &= \ \varphi_1 + (\theta^*_{r+d}-\theta^*_r)(\theta_{t+1} - \theta_t),\\
\varphi_d \ &= \ \phi_1 + (\theta^*_{r+d}-\theta^*_r)(\theta_{t+d-1} - \theta_{t+d}).
\end{align*} 
\end{lemma}

\proof Immediate from (\ref{varphi}) and (\ref{phi}). \hfill $\Box$



\begin{lemma}\label{previous}
Suppose that $W$ and $W'$ are thin irreducible $T\!$-modules with the same endpoint, dual endpoint and diameter $d >0$.  Then the following are equivalent:
\begin{align*}
 \varphi_1(W) &=  \varphi_1(W'), & \varphi_d(W) &= \varphi_d(W'),\\
 \phi_1(W) &= \phi_1(W'), &  \phi_d(W) &= \phi_d(W'),\\
a_0(W) &= a_0(W'), &  a_d(W) &= a_d(W'),\\
 a^{*}_0(W) &= a^{*}_0(W'), & a^{*}_d(W) &= a^{*}_d(W').
\end{align*}
\end{lemma}

\proof
Combine Lemma \ref{mgaphi}, (\ref{a01}), (\ref{a02}), (\ref{a01*}), (\ref{a02*}).   \hfill $\Box$

\begin{theorem}
Suppose that $W$ and $W'$ are thin irreducible $T\!$-modules with common diameter $d$.  
\begin{enumerate}
\item[\rm(i)] Assume $d=0$.  Then $W$ and $W'$ are isomorphic as $T\!$-modules if and only if they have the same endpoint and dual endpoint. 
\item[\rm(ii)] Assume $d>0$.  Then $W$ and $W'$ are isomorphic as $T\!$-modules if and only if they have the same endpoint, dual endpoint and all of the quantities in Lemma \ref{previous}.
\end{enumerate}
\end{theorem}

\proof
(i)  Immediate from Lemma \ref{isoclass}.\\
(ii) By Theorem \ref{formofintarray2}, $W$ and $W'$ have the same intersection numbers (resp. dual intersection numbers) if and only if they have the same endpoint, dual endpoint, diameter and $\varphi_1(W) = \varphi_1(W')$. Combining this with Lemmas \ref{isoclass}, \ref{previous}, we obtain the desired result.   \hfill $\Box$      

\section{Two examples of $Q$-polynomial distance-regular graphs}
In this section, we apply the results that we have obtained in Section \ref{everything} to several examples of $Q$-polynomial distance-regular graphs.  We will continue talking about the $T$-module $W$ in Assumption \ref{W} but now we will impose extra conditions on $\Gamma$.

\begin{definition} \rm\label{typeI}
The graph $\Gamma$ is said to have \textit{$q$-Racah type} whenever its parameter array $(\{\theta_i\}_{i=0}^D, \{\theta^*_i\}_{i=0}^D,$ $\{\varphi_i \}_{i=1}^D, \{\phi_i \}_{i=1}^D)$ satisfy the following.\\
For $0 \leq i \leq D$,
\begin{align*}
\theta_{i} \ &= \ \theta_0 + hq^{-i}(1-q^i)(1-sq^{i+1}), \\
\theta^*_{i} \ &= \ \theta^*_0 + h^*q^{-i}(1-q^i)(1-s^*q^{i+1}).
\end{align*}
For $1 \leq i \leq D$,
\begin{align*}
\varphi_i \ & = \ hh^*q^{1-2i}(1-q^i)(1-q^{i-D-1})(1-r_1q^i)(1-r_2q^i),\\
\phi_i \ &= \ hh^*q^{1-2i}(1-q^i)(1-q^{i-D-1})(r_1-s^*q^i)(r_2 -s^*q^i)/s^*.
\end{align*}
In the above, $q,h, h^*, r_1, r_2, s, s^*$ are complex scalars such that $r_1r_2 = ss^*q^{D+1}$, $hh^*ss^* \neq 0, \ q \notin \{-1, 0, 1 \}$.  
\end{definition}

\begin{lemma} \label{thetadiff}
Let $\Gamma$ be as in Definition \ref{typeI}.  Then for $0 \leq i, j \leq D$,
\begin{align*}
\theta_{i} - \theta_{j} \ &= \ h(q^i - q^j)(sq - q^{-i-j}), \\
\theta^*_{i} - \theta^*_{j} \ &= \ h^*(q^i - q^j)(s^*q - q^{-i-j}). 
\end{align*}
\end{lemma}

\proof
Routine calculation using Definition \ref{typeI}. \hfill $\Box$

\begin{lemma}
With reference to Definition \ref{typeI}, none of $q^i, r_1q^i, r_2q^i, s^*q^i/r_1, s^*q^i/r_2$ is equal to $1$ for $1 \leq i \leq D$.  Moreover, neither of $sq^i, s^*q^i$ is equal to $1$ for $2 \leq i \leq 2D$.  
\end{lemma}

\proof
The first assertion follows from Definition \ref{typeI} and the fact that for $1 \leq i \leq D$, $\theta_i \neq \theta_0$, $\theta^*_i \neq \theta^*_0$, $\varphi_i \neq 0, \ \phi_i \neq 0$.  The second assertion is immediate from Lemma \ref{thetadiff} and the fact that the eigenvalues (resp. dual eigenvalues) of $\Gamma$ are mutually distinct. \hfill $\Box$

\begin{lemma}\label{lemma:tauW}
Let $\Gamma$ be as in Definition \ref{typeI}.  Let $\{\varphi_i(W)\}_{i=1}^d$ and $\{\phi_i(W)\}_{i=1}^d$ be the first split sequence and second split  sequence of $W$\!, respectively.  Then there exists $\tau(W) \in \mathbb{C}$ such that for $1 \leq i \leq d$,
\begin{align}
\varphi_i(W) \ &= \ hh^*(1-q^{i})(1-q^{d-i+1})(\tau(W)- ss^*q^{r+t+i+1} - q^{-r-t-i-d}), \label{varphifirst}\\
\phi_i(W) \ &= \ hh^*(1-q^{i})(1-q^{d-i+1})(\tau(W)- s^*q^{r-t-d+i}- sq^{t-r-i+1}). \label{phifirst}
\end{align}   
\end{lemma}

\proof
Since $h,  h^*$ are both nonzero and $q, q^d$ are both not equal to $1$, there exists $\tau(W)$ such that (\ref{varphifirst}) holds for $i=1$.  Plugging $\varphi_1(W)$ in (\ref{phi}) and using Lemma \ref{thetadiff}, we routinely obtain that (\ref{phifirst}) holds for $1 \leq i \leq d$. Evaluating (\ref{varphi}) using (\ref{phifirst}) with $i=1$ and repeating the same argument above, we find that (\ref{varphifirst}) holds for $1 \leq i \leq d$. \hfill $\Box$\\

\indent
We make a comment about our notation used in Lemma \ref{lemma:tauW}. In the proof of \cite[Theorem 35.15]{madrid}, there are scalars $\tau, h, h^*$.  Our present $h, h^*$ are the same as those in \cite[Theorem 35.15]{madrid}.  However, our $\tau(W)$ is equal to $\tau/hh^*$.     

\begin{theorem}\label{splitfortypeI}
Let $\Gamma$ be as in Definition \ref{typeI}.  Let $r_1(W), r_2(W)$ be the roots of 
$$\lambda^2 - \tau(W)q^{r+t+d}\lambda + ss^*q^{2r+2t+d+1} =0,$$
where $\tau(W)$ is from Lemma \ref{lemma:tauW}.
Then for $1 \leq i \leq d$,
\begin{align}
\varphi_i(W) \ &= \ hh^*q^{1-2i-t-r}(1-q^i)(1-q^{i-d-1})(1-r_1(W)q^i)(1-r_2(W)q^i), \label{typeIvarphi}\\
\phi_i(W) \ &= \ hh^*q^{1-2i-t-r}(1-q^i)(1-q^{i-d-1})(r_1(W)-s^*q^{i+2r})(r_2(W) -s^*q^{i+2r})/s^*q^{2r}\!. \label{typeIphi}
\end{align}
\end{theorem}

\proof 
Note that 
\begin{equation}
r_1(W)r_2(W) = ss^*q^{2r+2t+d+1}, \qquad r_1(W) + r_2(W) = \tau(W)q^{r+t+d}. \label{fortau(W)}
\end{equation}
Eliminating $\tau(W), ss^*$ in (\ref{varphifirst}) using (\ref{fortau(W)}), we routinely obtain (\ref{typeIvarphi}).  Arguing similarly, we obtain (\ref{typeIphi}).\hfill $\Box$\\

\indent
The next theorem will involve basic hypergeometric series.  For the definition, see \cite[p.4]{Gasper}.

\begin{theorem}\label{uitypeI}
Let $\Gamma$ be as in Definition \ref{typeI}.  Then
$$u_i(\theta_{t+j}) \ = \ {}_4\phi_3 \left( \!\begin{array}{c} q^{-i}, s^*q^{2r+i+1}, q^{-j}, sq^{2t+j+1}\\r_1(W)q, r_2(W)q, q^{-d} \end{array} \Bigg \vert \; q, q \right) \qquad (0 \leq i, j \leq d), $$
where $u_i = u_i^W, \ r_1(W),\ r_2(W)$ are from Definitions \ref{defdefui}, \ref{splitfortypeI}.
\end{theorem}

\proof Routine calculation using (\ref{uitaui}) and Lemmas \ref{thetadiff}, \ref{splitfortypeI}.\hfill $\Box$\\

\indent
The polynomials $u_i$ are $q$-Racah polynomials.  For the definition of $q$-Racah polynomials, see \cite{orthogpoly}.   

\begin{theorem}\label{intypeI}
Let $\Gamma$ be as in Definition \ref{typeI}. Then the intersection numbers of $W$ are as follows:  
\begin{align*}
b_0(W) \ &= \ \frac{hq^{-t}(1-q^{-d})(1-r_1(W)q)(1-r_2(W)q)}{(1-s^*q^{2r+2})} \\
b_i(W) \ &= \ \frac{hq^{-t}(1-q^{i-d})(1-s^*q^{2r+i+1})(1-r_1(W)q^{i+1})(1-r_2(W)q^{i+1})}{(1-s^*q^{2r+2i+1})(1-s^*q^{2r+2i+2})} \qquad \ \ (1 \leq i \leq d-1),\\
c_i(W) \ & = \ \frac{hq^{-t}(1-q^i)(1-s^*q^{2r+i+d+1})(r_1(W) - s^*q^{2r+i})(r_2(W) - s^*q^{2r+i})}{s^*q^{2r+d}(1-s^*q^{2r+2i})(1-s^*q^{2r+2i+1})} \quad (1 \leq i \leq d-1),\\
c_d(W) \ &= \ \frac{hq^{-t}(1-q^d)(r_1(W) - s^*q^{2r+d})(r_2(W) - s^*q^{2r+d})}{s^*q^{2r+d}(1-s^*q^{2r+2d})}\\
a_i(W)\ &= \ \theta_t - b_i(W) - c_i(W) \qquad \qquad \qquad \qquad \qquad \qquad\qquad \qquad\qquad \qquad \quad \  (0 \leq i \leq d),
\end{align*}
where $r_1(W), r_2(W)$ are from Theorem \ref{splitfortypeI}.  To obtain the dual intersection numbers of $W$, replace $(h, s^*\!, r, t)$ with $(h^*\!, s, t, r)$.
\end{theorem}

\proof
Evaluate the equations on the left in (\ref{biprod}) and (\ref{ciprod}) using  Lemma \ref{thetadiff} and Theorem \ref{splitfortypeI}. \hfill $\Box$

\begin{corollary}
Let $\Gamma$ be as in Definition \ref{typeI}. Then the intersection numbers of $\Gamma$ are as follows:  
\begin{align*}
b_0 \ &= \ \frac{h(1-q^{-D})(1-r_1q)(1-r_2q)}{(1-s^*q^{2})} \\
b_i \ &= \ \frac{h(1-q^{i-D})(1-s^*q^{i+1})(1-r_1q^{i+1})(1-r_2q^{i+1})}{(1-s^*q^{2i+1})(1-s^*q^{2i+2})} \qquad \ \ (1 \leq i \leq D-1),\\
c_i \ & = \ \frac{h(1-q^i)(1-s^*q^{i+D+1})(r_1 - s^*q^{i})(r_2 - s^*q^{i})}{s^*q^{D}(1-s^*q^{2i})(1-s^*q^{2i+1})} \qquad \quad \ \ (1 \leq i \leq D-1),\\
c_D \ &= \ \frac{h(1-q^D)(r_1 - s^*q^{D})(r_2 - s^*q^{D})}{s^*q^{D}(1-s^*q^{2D})}\\
a_i \ &= \ \theta_0 - b_i - c_i \qquad \qquad \qquad \qquad \qquad \qquad\quad \ \qquad \qquad \quad \  (0 \leq i \leq D),
\end{align*}
where $r_1$, $r_2$ are from Definition \ref{typeI}.
To obtain the dual intersection numbers of $\Gamma$, replace $(h, s^*)$ with $(h^*, s)$.
\end{corollary}

\proof
Apply Theorem \ref{intypeI} with $W$ equal to the trivial $T$-module and use Lemma \ref{fortrivial}. \hfill $\Box$\\

We now turn our attention to graphs with classical parameters.  

\begin{definition} \rm\label{classical}
Let $b, \alpha, \sigma \in \mathbb{C}$ with $b \notin \{-1,0,1\}$.  The graph $\Gamma$ is said to have \textit{classical parameters} $(D, b, \alpha, \sigma)$ whenever 
\begin{align*}
c_i \ &= \ \frac{b^i - 1}{b-1}\left(1 + \alpha \frac{b^{i-1} -1}{b-1}  \right) \qquad (1 \leq i \leq D),\\
b_i \ &= \ \frac{b^D-b^i}{b-1}\left( \sigma - \alpha  \frac{b^i-1}{b-1}  \right) \qquad (0 \leq i \leq D-1).
\end{align*}
\end{definition}

\begin{theorem}{\rm\cite[Corollary 8.4.4]{BCN}}\label{resultclassical}
Let $\Gamma$ be as in Definition \ref{classical}. The following hold.
\begin{enumerate}
\item[\rm(i)]  There exists an ordering $\{\theta_i \}_{i=0}^D$ of the eigenvalues such that for $0 \leq i \leq D$,  
$$\theta_i = \eta + \mu b^i + h b^{-i}\mspace{-1mu},$$
where 
\begin{align*}
\eta \ &= \ \frac{(\sigma -1)(1-b) - \alpha(b^D+1)}{(b-1)^2},\\
\mu \ &=\ \frac{\alpha - b + 1}{(b-1)^2},\\
h \ &= \ \frac{b^D(\sigma b - \sigma + \alpha)}{(b-1)^2}.
\end{align*}
\item[\rm(ii)]  $\Gamma$ is $Q$-polynomial with respect to $\{ \theta_i \}_{i=0}^D$. 
\end{enumerate}
\end{theorem}

\indent  
Let $E$ be the primitive idempotent of $\Gamma$ corresponding to $\theta_1$.  Let $\{\theta^*_i \}_{i=0}^D$ be its corresponding dual eigenvalue sequence.  Our next goal is to express these values in terms of $\alpha, b, \sigma, D$.  

\begin{theorem}{\rm\cite[Corollary 8.4.4]{BCN}} \label{lemma:deval1}
Let $\{\theta^*_i \}_{i=0}^D$ be the dual eigenvalues corresponding to $E$.  Then for $0 \leq i \leq D$,
\begin{equation}
\theta^*_i = \eta^* + h^*b^{-i}\mspace{-1mu}, \label{dualclass}
\end{equation}
where 
\begin{align}
\eta^* \  &= \ \theta^*_0\left(1 + \frac{b}{b-1}\frac{(\sigma - \alpha)(b^{D-1}-1) - b + 1- \sigma (b^D-1)}{\sigma(b^D-1)}\right), \nonumber\\
h^* \  &= \ \theta^*_0 - \eta^*\!. \label{eq:hstar} 
\end{align}
\end{theorem}

\indent
Observe that by (\ref{dualclass}), $h^* \neq 0$ since $\theta^*_i \neq \theta^*_0$ for $1 \leq i \leq D$.  Later in this section, we will express $\theta^*_0$ in terms of $\alpha, b, \sigma, D$.

\begin{lemma}\label{diffclassical}
Let $\Gamma$ be as in Definition \ref{classical}.  Then for $0 \leq i, j \leq D$,
\begin{align*}
\theta_{i} - \theta_{j} \  &= \ (b^i - b^j)(\mu - hb^{-i-j}), \\
\theta^*_{i} - \theta^*_{j} \  &= \ h^*b^{-i-j}(b^j - b^i), 
\end{align*}
where $\mu, \ h, \ h^*$ are the from Theorems \ref{resultclassical}, \ref{lemma:deval1}. 

\end{lemma}

\proof Routine calculation using Theorems \ref{resultclassical}, \ref{lemma:deval1}. \hfill $\Box$

\begin{lemma} \label{splitclass}
Let $\Gamma$ be as in Definition \ref{classical}. Then there exists $\tau(W) \in \mathbb{C}$ such that the first split sequence and second split sequence of $W$ are given by
\begin{align}
\varphi_i(W) \  &= \ (1-b^i)(1-b^{d-i+1})(\tau(W) - hh^* b^{-r-t-i-d}) \qquad  (1 \leq i \leq d) \label{varphiclass},\\
\phi_i(W) \  &= \ (1-b^i)(1-b^{d-i+1})(\tau(W) - h^*\mu b^{-r+t-i}) \qquad \quad \  (1 \leq i \leq d)\label{phiclass},
\end{align}
where $\mu, \ h, \ h^*$ are from  Theorems \ref{resultclassical}, \ref{lemma:deval1}.
\end{lemma}

\proof
Similar to the proof of Lemma \ref{lemma:tauW}.\hfill $\Box$\\

\indent
Applying Lemma \ref{splitclass} with $W$ equal to the trivial $T$-module and using Definition \ref{PAgamma}, we obtain the following corollary.  

\begin{corollary}\label{cor:splitclass}
Let $\Gamma$ be as in Definition \ref{classical}. Then the first split sequence and second split sequence of $\Gamma$ are as follows: 
\begin{align*}
\varphi_i \  &= \ (1-b^i)(1-b^{D-i+1})(\tau - hh^* b^{-i-D}) \qquad (1 \leq i \leq D), \\
\phi_i \  &= \ (1-b^i)(1-b^{D-i+1})(\tau - h^*\mu b^{-i}) \qquad \ \ \ (1 \leq i \leq D),
\end{align*} 
where $\mu,h, h^*$ are from Theorems \ref{resultclassical}, \ref{lemma:deval1} and $\tau$ is the $\tau(W)$ associated with the trivial $T$-module. 
\end{corollary}

\indent 
Observe that the parameter $h$ given in Theorem \ref{resultclassical} may or may not be zero.  Consider a thin irreducible $T$-module $W\!.$  Note that if $h = 0$, then $\tau(W) \neq 0$.  This follows from (\ref{varphiclass}) and the fact that $\varphi_i(W) \neq 0$ for $1 \leq i \leq d$.

\begin{theorem}
Let $\Gamma$ be as in Definition \ref{classical}.  For $0 \leq i, j \leq d$,
$$
u_i(\theta_{t+j}) \  = \ \left\{
\begin{array}{lc}
{}_3\phi_2 {\displaystyle\Biggl({{b^{-i}, \;b^{-j},\;\frac{\mu b^{2t+j}}{h}}\atop
{b^{-d},\;\; \frac{\tau(W)b^{r+t+d+1}}{hh^*}}}\;\Bigg\vert \; b,\;b\Biggr)}  & \mbox{ if }h \neq 0\\
{}_2\phi_1 {\displaystyle \Biggl({{b^{-i}, \;b^{-j}}\atop
{b^{-d}}}\;\Bigg\vert \; b,\;\frac{\mu h^*b^{-r+t+j-d}}{\tau(W)}\Biggr)} & \mbox{ if }h = 0,
\end{array}
\right.
$$
where $u_i = u_i^W$ is from Definition \ref{defdefui} and $\mu,h, h^*, \tau(W)$ are from Theorems \ref{resultclassical}, \ref{lemma:deval1} and Lemma \ref{splitclass}.
\end{theorem}

\proof Routine calculation using (\ref{uitaui}) and Lemmas \ref{diffclassical}, \ref{splitclass}. \hfill $\Box$

\begin{theorem} \label{intclass}
Let $\Gamma$ be as in Definition \ref{classical}. Then the intersection numbers of $W$ are as follows:
\begin{align}
b_i(W) \  &= \ b^{r+2i+1}(1-b^{d-i})(\tau(W) - hh^*b^{-r-t-i-d-1})/h^* \qquad (0 \leq i \leq d-1), \label{eq:biclass}\\
c_i(W) \  &= \ b^{r+i}(b^i-1)(\tau(W) - h^*\mu b^{-r+t-i})/h^* \qquad \qquad \qquad \  (1 \leq i \leq d), \label{eq:ciclass}\\
a_i(W) \  &= \ \theta_t - b_i(W) - c_i(W) \qquad \qquad\qquad\qquad\qquad\qquad \quad \ \ (0 \leq i \leq d), \nonumber
\end{align}
where $\mu,h, h^*, \tau(W)$ are from Theorems \ref{resultclassical}, \ref{lemma:deval1} and Lemma \ref{splitclass}.
\end{theorem}

\proof
Evaluate the equations on the left in (\ref{biprod}) and (\ref{ciprod}) using Lemmas \ref{diffclassical}, \ref{splitclass}.  \hfill $\Box$

\begin{theorem} \label{thm:starclass}
Let $\Gamma$ be as in Definition \ref{classical}. Then the dual intersection numbers of $W$ are as follows:
\begin{align}
b^*_0(W) \  &= \ \frac{(b^{d}-1)(\tau(W) - hh^*b^{-r-t-d-1})}{\mu b^t - hb^{-t-1}}, \nonumber\\
b^*_i(W) \  &= \ \frac{b^{-i}(b^{d-i}-1)(\tau(W) - hh^*b^{-r-t-i-d-1})(\mu b^t - hb^{-t-i})}{(\mu b^t - hb^{-t-2i-1})(\mu b^t - hb^{-t-2i})} \quad \quad \quad (1 \leq i \leq d-1), \label{bi*W}\\
c^*_i(W) \  &= \ \frac{b^{d-2i+1}(1-b^{i})(\tau(W) - h^*\mu b^{-r+t-d+i-1})(\mu b^{t}- hb^{-t-i-d})}{(\mu b^{t}- hb^{-t-2i})(\mu b^{t}- hb^{-t-2i+1})} \quad (1 \leq i \leq d-1 ),\label{ci*W} \\
c^*_d(W) \  &= \ \frac{b^{-d+1}(1-b^{d})(\tau(W) - h^*\mu b^{-r+t-1})}{\mu b^{t}- hb^{-t-2d+1}},\nonumber \\
a^*_i(W) \  &= \ \theta^*_r - b^*_i(W) - c^*_i(W) \qquad \qquad \qquad \qquad \qquad \qquad \quad\quad\quad\quad\quad\quad (0 \leq i \leq d),\nonumber
\end{align}
whete $\mu,h, h^*, \tau(W)$ are from Theorems \ref{resultclassical}, \ref{lemma:deval1} and Lemma \ref{splitclass}.
\end{theorem}

\proof
Note that for $1 \leq i \leq d-1$, $\mu b^t - hb^{-t-2i-1} \neq 0$ since this is a factor of $\theta_{t+i} - \theta_{t+i+1}$ by Lemma \ref{diffclassical} and the eigenvalues of $\Gamma$ are mutually distinct. Similarly, $\mu b^t - hb^{-t-2i}$ and $\mu b^t - hb^{-t-2i+1}$ are nonzero since these are factors of $\theta_{t+i-1} - \theta_{t+i+1}$ and $\theta_{t+i} - \theta_{t+i-1}$, respectively.  Arguing as in the proof of Theorem \ref{intclass}, we obtain the desired result. \hfill $\Box$\\

\indent
In Definition \ref{classical}, we gave a formula for the intersection numbers of $\Gamma$ in terms of $\alpha, b, \sigma, D$.  We now give an alternate formula in terms of $\mu, h, h^*$.

\begin{theorem}
Let $\Gamma$ be as in Definition \ref{classical}.  Then the intersection numbers of $\Gamma$ are as follows:
\begin{align}
b_i \  &= \ b^{2i+1}(1-b^{D-i})(\tau - hh^*b^{-i-D-1})/h^*  \qquad (0 \leq i \leq D-1),\nonumber\\
c_i\  &= \ b^{i}(b^i-1)(\tau - h^*\mu b^{-i})/h^* \label{eq:citrivial} \qquad \qquad\qquad \quad (1 \leq i \leq D),\\
a_i \  &= \ \theta_0 - b_i - c_i \qquad \qquad\qquad\qquad\qquad\quad\qquad   \ (0 \leq i \leq D), \nonumber
\end{align}
where $\mu,h, h^*, \tau$ are from Theorems \ref{resultclassical}, \ref{lemma:deval1} and Corollary \ref{cor:splitclass}. 
\end{theorem}

\proof Immediate from Lemmas \ref{fortrivial} and \ref{intclass}. \hfill $\Box$\\

\indent
We now give a formula for the dual intersection numbers of $\Gamma$.
\begin{theorem}
Let $\Gamma$ be as in Definition \ref{classical}. Then the dual intersection numbers of $\Gamma$ are as follows:
\begin{align}
b^*_0 \  &= \ \frac{(b^{D}-1)(\tau - hh^*b^{-D-1})}{\mu  - hb^{-1}}, \nonumber\\
b^*_i \  &= \ \frac{b^{-i}(b^{D-i}-1)(\tau - hh^*b^{-i-D-1})(\mu  - hb^{-i})}{(\mu - hb^{-2i-1})(\mu  - hb^{-2i})} \qquad \qquad \qquad (1 \leq i \leq D-1), \nonumber\\
c^*_i \  &= \ \frac{b^{D-2i+1}(1-b^{i})(\tau - h^*\mu b^{-D+i-1})(\mu - hb^{-i-D})}{(\mu - hb^{-2i})(\mu - hb^{-2i+1})} \qquad \qquad (1 \leq i \leq D-1),\label{eq:cistartrivial}\\
c^*_D \  &= \ \frac{b^{-D+1}(1-b^{D})(\tau - h^*\mu b^{-1})}{(\mu - hb^{-2D+1})},\nonumber\\
a^*_i \  &= \ \theta^*_0 - b^*_i - c^*_i  \ \qquad \quad\qquad \qquad\qquad \qquad\qquad \qquad\qquad \qquad (0 \leq i \leq D),\nonumber
\end{align}
where $\mu,h, h^*, \tau$ are from Theorems \ref{resultclassical}, \ref{lemma:deval1} and Corollary \ref{cor:splitclass}. 
\end{theorem}

\proof
Immediate from Lemma \ref{fortrivial} and Theorem \ref{thm:starclass} \hfill $\Box$\\

\indent
Recall that in Lemma \ref{lemma:deval1}, we gave a formula for the dual eigenvalues of $\Gamma$ in terms of $\alpha, b, \sigma, \theta^*_0$.  We are now ready to solve for $\theta^*_0$.

\begin{lemma}\label{lemma:theta0*}
Let $\Gamma$ be as in Definition \ref{classical}.  Let $h^*$ and $\tau$ be as in Theorem \ref{lemma:deval1} and Corollary \ref{cor:splitclass}, respectively.  Then
\begin{align}
h^* \  &= \ -\frac{b^{D-1}(\mu b^2 - h)(\mu b - h)}{(\mu b^{D+1} - h)(b^{D-1} + \mu( b -1)(b^{D-1}-1))},\label{eq:h*}\\
\tau \  &= \ \frac{h^*(1+ \mu b- \mu)}{b(b-1)}, \label{eq:tau}\\
\theta^*_0 \  &= \ \frac{h^*\sigma(b^D-1)(1-b)}{b((\sigma - \alpha)(b^{D-1}-1) - b+1 - \sigma(b^D-1))}. \label{eq:theta0*}
\end{align} 
\end{lemma}

\proof
To obtain (\ref{eq:h*}) and (\ref{eq:tau}), solve the system of equations in (\ref{eq:citrivial}) and (\ref{eq:cistartrivial}) with  $i=1$ and use the fact that $c_1 =1=c^*_1$. Line (\ref{eq:theta0*}) is immediate from (\ref{eq:hstar}). \hfill $\Box$     \\

\indent
In Theorem \ref{intclass}, we gave a formula for the intersection numbers of $W\!.$  We now give an alternate formula which is reminiscent of Definition \ref{classical}
\begin{theorem}
Let $\Gamma$ be as in Definition \ref{classical}. Then 
\begin{align}
c_i(W) \  &= \ \frac{b^i - 1}{b-1}\left(c_1(W) + \alpha(W) \frac{b^{i-1} -1}{b-1}  \right) \quad (1 \leq i \leq d), \label{eq:ciclass2}\\
b_i(W) \  &= \ \frac{b^d-b^i}{b-1}\left( \sigma(W) - \alpha(W)  \frac{b^i-1}{b-1}  \right) \qquad (0 \leq i \leq d-1),\label{eq:biclass2}
\end{align}
where
\begin{align*}
\alpha(W) \  &= \ \tau(W)b^{r+1}(b-1)^2/h^*,\\
\sigma(W) \  &= \ \frac{hb^{-t}(b-1)^2 - \alpha(W)b^d}{b^d(b-1)}.
\end{align*}
\end{theorem}

\proof
Comparing the right side of (\ref{eq:ciclass2}) with that of (\ref{eq:ciclass}), we find that (\ref{eq:ciclass2}) holds. Comparing the right side of (\ref{eq:biclass2}) with that of (\ref{eq:biclass}), we obtain (\ref{eq:biclass2}). \hfill $\Box$    

\bigskip

\noindent{\Large\bf Acknowledgments}

\medskip
\noindent
This paper was written while the author was an honorary fellow at the University of Wisconsin-Madison, January-December 2009, with support from HEDP-FDP Sandwich Program of the Commission on Higher Education, Philippines.  The author is greatly indebted to Professor Terwilliger for his many valuable ideas and suggestions.

\bigskip
\noindent Diana R. Cerzo \hfil\break
\noindent Institute of Mathematics \hfil\break
\noindent University of the Philippines \hfil\break
\noindent C.P. Garcia St., Diliman \hfil\break
\noindent Quezon City, Philippines 1101 \hfil\break
\noindent email:  {\tt drcerzo@math.upd.edu.ph}

\end{document}